\documentclass[11pt]{amsart}
\usepackage[T1]{fontenc}
\usepackage{lmodern,amsmath,amssymb,amsthm,mathtools}
\usepackage[margin=1.08in]{geometry}
\usepackage{microtype}

\newtheorem{theorem}{Theorem}[section]
\newtheorem{lemma}[theorem]{Lemma}
\newtheorem{proposition}[theorem]{Proposition}
\newtheorem{corollary}[theorem]{Corollary}
\theoremstyle{definition}
\newtheorem{example}[theorem]{Example}
\newtheorem{remark}[theorem]{Remark}

\newcommand{\R}{\mathbb R}
\newcommand{\E}{\mathbb E}
\newcommand{\Hil}{\mathcal H}
\DeclareMathOperator{\tr}{tr}
\DeclareMathOperator{\Ker}{ker}

\title[Continuous Gaussian H\"older inequalities]{Continuous improved H\"older and reverse H\"older inequalities for Gaussian random vectors}
\author[A.-J. Li]{Ai-Jun Li}
\address[A.-J. Li]{School of Science,
	Zhejiang University of Science and Technology, Hangzhou 310023, China}
\email{liaijun72@163.com}	
\author[Q. Huang]{Qingzhong Huang}
\address[Q. Huang]{College of data science, Jiaxing University, Jiaxing 314001,	China} \email{hqz376560571@163.com} 
\begin{document}

\begin{abstract}
We establish continuous versions of the improved H\"older and reverse H\"older inequalities for correlated Gaussian random vectors, where a finite covariance matrix is replaced by a covariance operator acting on a Hilbert space of parameter fields. Under operator comparisons between the covariance and the exponent multiplier, we prove both inequalities and characterize equality in terms of a compatible probability density on the underlying Gaussian space. We further describe this density as a product of independent factors on common eigenspaces and a Gaussian factor whose mean and covariance satisfy explicit algebraic conditions. In the reverse setting, positive and negative exponents may coexist. The continuous product is defined by ordered truncation, allowing zeros in the positive-exponent factors. The equality proof combines an exact Ornstein--Uhlenbeck semigroup comparison, a spectral gap away from the kernel of the governing operator difference, and characteristic-function deconvolution. We also treat exponential quadratic functions without requiring a uniform positive-definiteness margin on the precision matrices.
\end{abstract}
\thanks{The first author was supported by NSFC (Grant No. 12571146 and 12231006). The second author was supported by  NSFC (Grant No. 11701219) and the Qin Shen Scholar Program of Jiaxing University.}
\maketitle

\section{Introduction}

H\"older's inequality bounds the expectation of a product of random variables solely in terms of their marginal moments: if $X_1, X_2$ are random vectors and $p_1, p_2 > 1$ satisfy $p_1^{-1} + p_2^{-1} = 1$, then
\begin{equation}\label{eq:intro-holder}
\E[f_1(X_1)f_2(X_2)] \le (\E f_1(X_1)^{p_1})^{1/p_1}(\E f_2(X_2)^{p_2})^{1/p_2}
\end{equation}
for all nonnegative measurable functions $f_1, f_2$. Because \eqref{eq:intro-holder} ignores the dependence structure between $X_1$ and $X_2$, the conjugacy constraint on $(p_1, p_2)$ can be relaxed substantially when $(X_1, X_2)$ is jointly Gaussian. For instance, {if $(X_1,X_2)$ is a jointly Gaussian vector whose components are standard normal random variables with correlation $\rho \in (-1,1)$}, then \eqref{eq:intro-holder} holds for any $p_1, p_2 \ge 1$ satisfying
\begin{equation}\label{eq:intro-twovariable}
(p_1 - 1)(p_2 - 1) \ge \rho^2,
\end{equation}
a condition reflecting the positive semidefiniteness of the matrix
\[
\begin{pmatrix} p_1 - 1 & -\rho \\ -\rho & p_2 - 1 \end{pmatrix}.
\]
{If the restriction $p_1,p_2\ge1$ is dropped, $p_1p_2\ne0$, and the displayed matrix is negative semidefinite, the corresponding lower bound holds for strictly positive functions with finite nonzero moments. The admissible exponents may include negative values.}

In the multivariate setting, Chen, Dafnis, and Paouris \cite{CDP} systematically investigated these covariance-dependent estimates. Let $X = (X_1, \dots, X_m)$ be a centered jointly Gaussian vector with $X_i \in \R^{n_i}$ and block covariance matrix $T = (T_{ij})_{1 \le i,j \le m}$. For nonzero exponents $p_1, \dots, p_m$, set $P = \operatorname{diag}(p_1 T_{11}, \dots, p_m T_{mm})$. Under appropriate integrability conditions, they proved that
\begin{equation}\label{eq:intro-cdp}
\E\prod_{i=1}^m f_i(X_i)
\begin{cases}
\displaystyle \le \prod_{i=1}^m (\E f_i(X_i)^{p_i})^{1/p_i}, & T \le P, \\[2mm]
\displaystyle \ge \prod_{i=1}^m (\E f_i(X_i)^{p_i})^{1/p_i}, & T \ge P,
\end{cases}
\end{equation}
and showed by testing on exponential functions that the operator inequalities $T \le P$ and $T \ge P$ are necessary for \eqref{eq:intro-cdp} to hold across all admissible inputs. They further linked the upper estimate in \eqref{eq:intro-cdp} to the geometric Brascamp--Lieb inequality for Gaussian measures, deducing applications to Gaussian hypercontractivity, reverse hypercontractivity, and sharp Young-type inequalities.

{The finite-index result of \cite{CDP} permits blocks of different dimensions and an arbitrary joint Gaussian covariance. In contrast, our continuous family is obtained from one finite-dimensional Gaussian vector through coisometries with a common output dimension. Thus our hypotheses do not cover the full covariance setting of \cite{CDP}. The main points of the present formulation are the parameter-space operator conditions, the ordered extension of the product in the mixed-sign case, and the equality characterization under the stated hypotheses.}

These inequalities sit within the broader framework of Brascamp--Lieb and functional inequalities, which unify H\"older, Loomis--Whitney, and sharp Young convolution inequalities \cite{Beckner,BL}. Lieb \cite{Lieb} showed that Gaussian inputs determine the best constants in a wide class of Brascamp--Lieb inequalities, while Bennett, Carbery, Christ, and Tao \cite{BCCT} analyzed finiteness, extremizability, and structural properties of the underlying linear data via heat-flow monotonicity. See also Ball \cite{Ball} for geometric applications to volume ratios and reverse isoperimetric inequalities. On the reverse side, Barthe \cite{Barthe,BartheYoung} established a reverse Brascamp--Lieb inequality via optimal transport, where functions are coupled through a supremum over a linear constraint. Although the lower bound in \eqref{eq:intro-cdp} is also often termed an inverse or reverse Brascamp--Lieb inequality, it differs structurally in that it estimates the expectation of a direct product evaluated along correlated Gaussian projections. Heat-flow approaches to such inverse and correlation inequalities were developed by Barthe and Cordero-Erausquin \cite{BCE} and Barthe, Cordero-Erausquin, Ledoux, and Maurey \cite{BCLM}, building on classical semigroup estimates such as Nelson's hypercontractivity \cite{Nelson}, Gross's logarithmic Sobolev inequality \cite{Gross}, and Borell's reverse hypercontractivity \cite{Borell,MOS}. Subsequent developments on Gaussian kernels, regularized inequalities, entropic dualities, and equality characterizations include \cite{Ammari,BH,BW,BN,BKX,CE,CLL,CourtadeLiu,Gardner,LCCV,Milman,NT,Valdimarsson}.

A natural extension of the discrete multilinear theory replaces a finite collection of linear maps by a measurable family indexed by a measure space. Continuous Brascamp--Lieb inequalities were introduced by Barthe \cite{BartheContinuous} for measures on the sphere, and developed further by Cordero-Erausquin and Maurey \cite{CEM} via Borell's stochastic control approach. Li, Xi, and Zhang \cite{LXX} studied inequalities connected with Grassmannian cosine transforms. In \cite{LiHuang}, the authors studied continuous geometric Brascamp--Lieb inequalities on Grassmannians, including a reverse inequality expressed through a barycentric supremum. The lower bound considered here concerns a direct product under a Gaussian measure and has a different formulation. Proposition~2.4 of \cite{CEM} gives a Gaussian integral comparison under an integrated pointwise hypothesis and a contraction between two $L^2$ parameter spaces. Its stated hypotheses include continuity and a growth bound on one family. Section~5 of \cite{CEM} treats different fiber dimensions. Our formulation starts from the covariance operator associated with a finite-dimensional common Gaussian source, allows measurable normalized densities under the hypotheses below, and considers an ordered product in the mixed-sign case.

Specifically, we replace the finite index set in \eqref{eq:intro-cdp} by a finite measure space $(\Omega, \mu)$ and consider a family of random vectors
\[
X_s = B_s Z, \qquad B_s B_s^* = I_k,
\]
where $Z$ is a standard Gaussian vector in $\R^n$ and each $B_s : \R^n \to \R^k$ is a coisometry. Thus the common Gaussian source is finite-dimensional, and every fiber $X_s$ is a standard Gaussian vector in $\R^k$. The joint covariance of $\{X_s\}_{s \in \Omega}$ is encoded by the finite-rank operator $C = A A^*$ on $\Hil = L^2(\Omega, \mu; \R^k)$, where $(Az)(s) = B_s z$, while the measurable exponent function $p : \Omega \to \R$ acts as a multiplication operator $P$. The operator inequalities $P \ge C$ and $C \ge P$ serve as the exact continuous analogues of the matrix conditions in \eqref{eq:intro-cdp}, recovering the normalized finite-dimensional criteria when $\mu$ is counting measure on a finite set.

{The connection with \cite{CEM} is more precise in the forward case. If $p>0$ and $P\ge C$, then, with $d\theta=p^{-1}d\mu$, the map $z\mapsto(s\mapsto B_s z)$ is a contraction from $\R^n$ to $L^2(\theta;\R^k)$. Conversely, this contraction condition implies $P\ge C$. Applying \cite[Proposition~5.2]{CEM} to this map gives the bounded continuous forward inequality. Approximation then gives the full forward inequality stated below, including functions with zeros; see Proposition~\ref{prop:cem-reduction}. Thus we do not claim the forward inequality itself as a new consequence of a continuous covariance comparison.

When the parameter space has finite support, both inequalities, including the mixed-sign lower bound, follow from the finite Gaussian comparison in \cite{CDP} after rescaling the atomic masses and exponents; see Proposition~\ref{prop:finite-cdp}.}

{The general mixed-sign reverse inequality can also be obtained from the finite comparison in \cite{CDP} by approximation. Proposition~\ref{prop:nonatomic-cdp} gives the reduction: on each negative-exponent parameter cell, we average the maps and adjust the exponent to compensate for the resulting change in covariance. The ordered truncation in Section~5 then covers zeros and divergent logarithms. Thus the inequalities themselves follow from earlier comparison principles together with the approximation arguments described here. The recovery of equality requires a separate argument.}

Moving from a finite product to a continuous parameter integral $\exp\int_\Omega p(s)^{-1}\log g_s(B_s z)\,d\mu(s)$ creates an essential analytical subtlety: when $p$ changes sign, the positive and negative parts of the integrand may diverge simultaneously on sets of positive Gaussian measure, even when every individual moment is finite. In the reverse regime $C \ge P$, however, the finite-rank nature of $C$ forces the positive-exponent set $\{p > 0\}$ to consist of at most finitely many atoms of $\mu$, whereas negative exponents may be distributed over a nonatomic part. This structural dichotomy allows us to define an extended continuous product $F_-$ by truncating the positive-exponent factors from above and the negative-exponent factors from below in opposite directions. In particular, whenever a positive-exponent factor vanishes, $F_-$ is set to zero even if the reciprocal contribution from the negative-exponent part diverges. As shown in Section~5, reversing the order of these truncations can alter the value of the product.

Our main results (Theorems~\ref{thm:forward} and~\ref{thm:reverse}) establish both the forward and reverse inequalities under $P \ge C$ and $C \ge P$, respectively, and characterize equality through the existence of a compatible probability density $r$ on $(\R^n, \gamma_n)$. Compatibility requires that for $Y \sim r\gamma_n$, the projection $B_s Y$ has law $g_s \gamma_k$ and is independent of $(I - B_s^* B_s)Y$ for $\mu$-almost every $s$, while $Y$ splits into an arbitrary compatible component on the subspace $V = \Ker(I - A^* P^{-1} A)$ and an independent standard Gaussian component on $V^\perp$. Theorem~\ref{thm:explicit} gives a more explicit description: arbitrary independent density factors occur on specified intersections of $V$ with common eigenspaces of the projections, and the remaining factor is Gaussian with constrained mean and covariance. In the reverse setting, we also determine the zero set of any compatible density directly from the finitely many positive-exponent marginals, and construct explicit non-Gaussian extremizers that vanish on a half-space.

To prove the inequalities, we first establish a differential monotonicity identity along the Ornstein--Uhlenbeck semigroup for bounded functions bounded away from zero, and then pass to general inputs via carefully ordered truncations. The equality analysis requires a different mechanism, as the truncated approximations need not be extremizers themselves. Instead, we prove a pointwise comparison between the semigroup applied to the normalized product and the product formed from the smoothed marginals. When equality holds at $t=0$, this forces exact equality at every positive time $t>0$. Exploiting the finite rank of $C$, we obtain a spectral gap for $P - C$ (or $C - P$) on the orthogonal complement of its kernel. Passing to the limit along bounded approximations shows that the smoothed gradient field lies in the kernel for almost every $(t,z)$, which integrates to a common smooth potential splitting additively across each projection. Finally, characteristic-function deconvolution transfers both the marginal identities and the orthogonal independence from positive time $t>0$ back to $t=0$ without imposing extra moment assumptions.

The remainder of the paper is structured as follows. Section~2 introduces the setup and main theorems. Section~3 proves the operator facts and compares the inequalities with earlier Gaussian comparisons. Section~4 derives the semigroup identity, and Section~5 proves the inequalities for general functions. Sections~6 and~7 establish the equality criterion and its support consequences. Section~8 describes the corresponding Gaussian and independent factors. Finally, Section~9 treats exponential quadratic functions without assuming a uniform positive-definiteness bound on the fiber precision matrices.

\section{Notation and main results}

Let $\Omega$ be a standard Borel space equipped with the completion of a finite Borel measure $\mu$, and let $1 \le k \le n$. Assume that $s \mapsto B_s$ is measurable, where each $B_s : \R^n \to \R^k$ is a linear map satisfying the coisometry condition
\begin{equation}\label{eq:coiso}
B_s B_s^* = I_k \qquad \text{for } \mu\text{-almost every } s.
\end{equation}
We write $\gamma_d$ for the standard Gaussian probability measure on $\R^d$. For $Z \sim \gamma_n$, the random vectors $X_s = B_s Z$ each follow the law $\gamma_k$. Define the Hilbert space $\Hil = L^2(\Omega, \mu; \R^k)$ and the bounded linear operators $A : \R^n \to \Hil$ and $C : \Hil \to \Hil$ by
\[
(Az)(s) = B_s z, \qquad C = A A^*.
\]
Let $p : \Omega \to \R$ be a measurable exponent function satisfying
\begin{equation}\label{eq:pbound}
0 < p_* \le |p(s)| \le p^* < \infty \qquad \text{for } \mu\text{-almost every } s,
\end{equation}
and let $P$ denote the bounded multiplication operator $(Ph)(s) = p(s)h(s)$ on $\Hil$. Set
\begin{equation}\label{eq:K}
Q_s = B_s^* B_s, \qquad K = A^* P^{-1} A = \int_\Omega \frac{Q_s}{p(s)}\,d\mu(s), \qquad V = \Ker(I - K).
\end{equation}
By \eqref{eq:coiso}, each $Q_s$ is an orthogonal projection on $\R^n$.

Let $(f_s)_{s \in \Omega}$ be a jointly measurable family of functions on $\R^k$ such that $f_s$ is nonnegative and finite $\gamma_k$-almost everywhere when $p(s) > 0$, and strictly positive and finite $\gamma_k$-almost everywhere when $p(s) < 0$. We assume that
\begin{equation}\label{eq:masses}
0 < m_s = \int_{\R^k} f_s^{p(s)}\,d\gamma_k < \infty \quad \text{for } \mu\text{-a.e. } s, \qquad
\int_\Omega \left|\frac{\log m_s}{p(s)}\right| d\mu(s) < \infty,
\end{equation}
and define the normalized probability densities with respect to $\gamma_k$ by
\begin{equation}\label{eq:g}
g_s = \frac{f_s^{p(s)}}{m_s}, \qquad \int_{\R^k} g_s\,d\gamma_k = 1.
\end{equation}
All identities in \eqref{eq:g} are understood for $\mu$-almost every $s$. When a formula requires values at every parameter, we work on one common full-measure set where the hypotheses hold and assign harmless fixed values on its complement: take a fixed coisometry for $B_s$, set $p(s)=p_*$, and set $f_s=g_s=m_s=1$. These changes do not affect any operator, parameter integral, or almost-everywhere equality in the results below. In particular, every positive-mass atom is retained.

In the forward setting, assume $p > 0$ almost everywhere and $P \ge C$ in the sense of self-adjoint operators on $\Hil$. We define the normalized forward product by
\begin{equation}\label{eq:forward}
F_+(z) = \exp\left(\int_\Omega \frac{1}{p(s)} \log g_s(B_s z)\,d\mu(s)\right),
\end{equation}
with the convention $\log 0 = -\infty$ and $e^{-\infty} = 0$. Since $\log^+ x \le x$ and $\int g_s\,d\gamma_k = 1$, Fubini's theorem ensures that the positive part of the integrand in \eqref{eq:forward} is $\mu$-integrable for $\gamma_n$-almost every $z$, so $F_+(z) \in [0,\infty)$ is well defined $\gamma_n$-almost everywhere.

In the reverse setting, assume $C \ge P$. As shown in Lemma~\ref{lem:atoms}, the set $\{s \in \Omega : p(s) > 0\}$ consists of at most finitely many atoms of $\mu$, which we denote by $s_1, \dots, s_m$ (modulo $\mu$-null sets). Write $\alpha_i = \mu(\{s_i\})/p(s_i)$, $B_i = B_{s_i}$, and $g_i = g_{s_i}$, and define the finite positive measure $d\nu = -p^{-1}\,d\mu$ on $\{p < 0\}$. For each integer $N \ge 1$, set
\begin{equation}\label{eq:reverseN}
F_N(z) = \prod_{i=1}^m \min\{g_i(B_i z), N\}^{\alpha_i}
\exp\left(-\int_{\{p < 0\}} \log\max\{g_s(B_s z), N^{-1}\}\,d\nu(s)\right),
\end{equation}
and define the extended normalized reverse product by
\begin{equation}\label{eq:reverse}
F_-(z) = \lim_{N \to \infty} F_N(z),
\end{equation}
where an empty product is taken to be $1$. Notice that $F_N(z)$ is nondecreasing in $N$: as $N$ increases, each factor $\min\{g_i, N\}^{\alpha_i}$ increases while $-\log\max\{g_s, N^{-1}\}$ also increases, so $F_N(z)$ is nondecreasing in $N$ and {the limit in \eqref{eq:reverse} exists in $[0,\infty]$ for $\gamma_n$-almost every $z$. On the remaining Gaussian null set, we set $F_N=F_-=0$ for every $N$.} In particular, if $g_i(B_i z) = 0$ for some positive-exponent atom $s_i$, then $F_N(z) = 0$ for every $N$ and hence $F_-(z) = 0$, regardless of whether the integral over $\{p < 0\}$ diverges to $-\infty$. Moreover, since $\log\max\{g_s, N^{-1}\} \ge -\log N$ and
\[
\int_{\R^n} \int_{\{p < 0\}} \log^+ g_s(B_s z)\,d\nu(s)\,d\gamma_n(z)
\le \int_{\{p < 0\}} \int_{\R^k} g_s\,d\gamma_k\,d\nu(s) = \nu(\{p < 0\}) < \infty,
\]
each $F_N$ is finite $\gamma_n$-almost everywhere, avoiding any indeterminate $\infty - \infty$ expressions.

To formulate the equality cases, we say that a probability measure $\rho = r\gamma_n$ on $\R^n$ (with density $r \ge 0$) is \emph{compatible} with $(B_s, p, g_s)$ if a random vector $Y \sim \rho$ satisfies:
\begin{enumerate}
\item $B_s Y$ has distribution $g_s \gamma_k$ for $\mu$-almost every $s$;
\item $Q_s Y$ and $(I - Q_s)Y$ are independent for $\mu$-almost every $s$;
\item $P_V Y$ and $(I - P_V)Y$ are independent, and $(I - P_V)Y$ is standard Gaussian on $V^\perp$, where $P_V$ denotes the orthogonal projection onto $V$.
\end{enumerate}

\begin{theorem}[Forward inequality]\label{thm:forward}
Assume $p > 0$ $\mu$-almost everywhere, $P \ge C$, and \eqref{eq:masses}. Then
\[
\int_{\R^n} F_+\,d\gamma_n \le 1.
\]
Equality holds if and only if there exists a compatible probability density $r$ on $(\R^n, \gamma_n)$ (in which case $r = F_+$ $\gamma_n$-almost everywhere).
\end{theorem}

\begin{theorem}[Reverse inequality]\label{thm:reverse}
Assume $C \ge P$ and \eqref{eq:masses}. Then
\[
\int_{\R^n} F_-\,d\gamma_n \ge 1,
\]
where the left-hand side may take the value $+\infty$. Equality holds if and only if there exists a compatible probability density $r$ on $(\R^n, \gamma_n)$ (in which case $r = F_-$ $\gamma_n$-almost everywhere).
\end{theorem}

In unnormalized form, the common right-hand side of both inequalities is
\[
\exp\left(\int_\Omega \frac{\log m_s}{p(s)}\,d\mu(s)\right) \in (0,\infty).
\]
Whenever $\int_\Omega |p(s)^{-1} \log g_s(B_s z)|\,d\mu(s) < \infty$, the extended product $F_-(z)$ coincides with 
$$\exp\int_\Omega p(s)^{-1}\log g_s(B_s z)\,d\mu(s),$$ while \eqref{eq:reverse} specifies the boundary behavior when a positive-exponent factor vanishes simultaneously with a negative-exponent divergence.

\begin{corollary}\label{cor:nonatomic}
Under $C \ge P$, we have $p < 0$ $\mu$-almost everywhere on the nonatomic part of $\mu$. In particular, if $\mu$ is nonatomic, equality holds in Theorem~\ref{thm:reverse} if and only if $g_s = 1$ $\gamma_k$-almost everywhere for $\mu$-almost every $s$.
\end{corollary}

\section{Operator facts}

We begin by recording the spectral and atomic consequences of the operator comparisons $P \ge C$ and $C \ge P$.

\begin{lemma}\label{lem:atoms}
Assume $C \ge P$. Then the set $\{s \in \Omega : p(s) > 0\}$ consists of at most finitely many atoms of $\mu$. If an atom $\{s_i\}$ has mass $\mu_i = \mu(\{s_i\})$ and exponent $p_i = p(s_i) > 0$, then $p_i \le \mu_i$. Consequently, there are at most $\mu(\Omega)/p_*$ positive-exponent atoms.
\end{lemma}

\begin{proof}
Since $A : \R^n \to \Hil$ has finite-dimensional domain, the operator $C = AA^*$ has finite rank and is therefore compact. Suppose toward a contradiction that there exists a nonatomic measurable subset $E \subset \Omega$ with $\mu(E) > 0$ on which $p \ge \varepsilon > 0$. Because $E$ is nonatomic, it supports an orthonormal sequence $(h_j)_{j \ge 1}$ in $\Hil$. Compactness of $C$ implies $\langle C h_j, h_j \rangle_{\Hil} \to 0$ as $j \to \infty$, whereas $\langle P h_j, h_j \rangle_{\Hil} \ge \varepsilon$ for all $j$, contradicting $C \ge P$. Taking $E = \{p \ge 1/j\}$ for $j \ge 1$ shows that $p \le 0$ almost everywhere on the nonatomic part of $\mu$, and by \eqref{eq:pbound}, $p < 0$ almost everywhere there.

Next, let $\{s_i\}$ be an atom of mass $\mu_i > 0$ with $p_i > 0$. Choose a unit vector $u \in \R^k$ and define $h \in \Hil$ by $h(s) = \mu_i^{-1/2} u \mathbf{1}_{\{s_i\}}(s)$, so $\|h\|_{\Hil} = 1$. Since $(A^* h) = \mu_i^{1/2} B_i^* u$ and $B_i B_i^* = I_k$, we have
\[
\langle C h, h \rangle_{\Hil} = |A^* h|^2 = \mu_i |B_i^* u|^2 = \mu_i, \qquad \langle P h, h \rangle_{\Hil} = p_i.
\]
Hence $C \ge P$ implies $p_i \le \mu_i$. Combining this with $p_i \ge p_*$ yields at most $\mu(\Omega)/p_*$ positive-exponent atoms.
\end{proof}

\begin{lemma}\label{lem:gap}
Suppose either $p > 0$ and $P \ge C$, or $C \ge P$. Set $T = P - C$ in the first case and $T = C - P$ in the second. Then $T \ge 0$, $\Ker T$ is finite-dimensional, and there exists a constant $c > 0$ such that
\begin{equation}\label{eq:gap}
\langle T h, h \rangle_{\Hil} \ge c \operatorname{dist}(h, \Ker T)^2 \qquad \text{for all } h \in \Hil.
\end{equation}
Moreover,
\begin{equation}\label{eq:kernel}
\Ker T = \{s \mapsto p(s)^{-1} B_s v : v \in V\}.
\end{equation}
\end{lemma}

\begin{proof}
In the forward case, write $T = P^{1/2}(I - P^{-1/2} C P^{-1/2})P^{1/2}$. Since $P^{-1/2} C P^{-1/2}$ is self-adjoint and finite rank, its nonzero spectrum consists of finitely many eigenvalues. Thus $1$ is either an isolated eigenvalue of finite multiplicity or outside the spectrum of $P^{-1/2} C P^{-1/2}$. In either event, $I - P^{-1/2} C P^{-1/2}$ has finite-dimensional kernel and a positive spectral gap on $(\Ker(I - P^{-1/2} C P^{-1/2}))^\perp$. Because $P^{1/2}$ and $P^{-1/2}$ are bounded invertible operators, this yields both the finite dimensionality of $\Ker T$ and \eqref{eq:gap}.

In the reverse case, let $E_+ \subset \Hil$ be the finite-dimensional subspace of fields supported on the positive-exponent atoms $\{s_1, \dots, s_m\}$. Define the bounded self-adjoint operator $D$ on $\Hil$ by $D = p_* I$ on $E_+$ and $D = -P$ on $E_+^\perp$. By \eqref{eq:pbound}, $D \ge p_* I$, so $D^{1/2}$ is bounded and invertible. Furthermore, $P + D = (P + p_* I)P_{E_+}$ has finite rank, whence
\[
D^{-1/2} T D^{-1/2} = D^{-1/2}(C - P)D^{-1/2} = I + R,
\]
where $R = D^{-1/2}(C - P - D)D^{-1/2}$ is self-adjoint and finite rank. Since $T \ge 0$, we have $I + R \ge 0$, and the finite-dimensional spectral theorem applied to $R$ gives a positive spectral gap for $I + R$ on $(\Ker(I + R))^\perp$. Conjugating by $D^{1/2}$ transfers this spectral gap and the finite dimensionality of the kernel to $T$.

Finally, $Th = 0$ is equivalent to $Ph = Ch = A(A^* h)$. Setting $v = A^* h \in \R^n$, this holds if and only if $h = P^{-1} A v$ and $v = A^* P^{-1} A v = K v$, i.e., $v \in V$.
\end{proof}

We will also use the following elementary closure property: any finite-dimensional subspace $W \subset L^2(\Omega, \mu; \R^k)$ is closed under convergence in $\mu$-measure to a finite measurable function. Indeed, if $(w_j) \subset W$ converges in $\mu$-measure to a finite measurable function $w$, then $\sup_j \|w_j\|_2 < \infty$. Otherwise, normalizing a subsequence with $\|w_{j_\ell}\|_2 \to \infty$ would yield a subsequence converging in $W$ (and hence in $\mu$-measure) to a unit vector $\tilde{w} \in W$, contradicting $w_{j_\ell}/\|w_{j_\ell}\|_2 \to 0$ in $\mu$-measure. Passing to a further subsequence converging in $W$ shows that $w \in W$ almost everywhere.

{\begin{proposition}[Comparison with \cite{CEM}]\label{prop:cem-reduction}
Assume $p>0$ almost everywhere. Set $d\theta=p^{-1}d\mu$ and define
\[
 (\mathcal A z)(s)=B_s z,\qquad
 \mathcal A:\R^n\longrightarrow L^2(\theta;\R^k).
\]
Then $P\ge C$ if and only if $\|\mathcal A\|\le1$. Under this condition, the inequality in Theorem~\ref{thm:forward}, including the case in which the densities $g_s$ vanish, follows from \cite[Proposition~5.2]{CEM} and the truncation argument in Section~5.
\end{proposition}

\begin{proof}
Since $P$ is positive and boundedly invertible,
\[
 P\ge AA^*
 \quad\Longleftrightarrow\quad
 P^{-1/2}AA^*P^{-1/2}\le I
 \quad\Longleftrightarrow\quad
 A^*P^{-1}A\le I.
\]
The last operator is $\mathcal A^*\mathcal A$, proving the equivalence.

We first take a uniformly bounded jointly measurable family $\psi_s$ whose members are continuous on $\R^k$. Put
\[
 \Psi(z)=\int_\Omega\psi_s(B_s z)\,d\theta(s).
\]
Then $\Psi$ is bounded and continuous. In the notation of \cite[Setting~2 and Proposition~5.2]{CEM}, take the first parameter space to be a single point of unit mass and the second to be $(\Omega,\theta)$, and take the input and output dimensions to be $n$ and $k$, respectively. The operator there is $\mathcal A$. It is a contraction and sends the constant vector $z$ to the field $s\mapsto B_s z$, as required by condition $(\mathrm{ii})_2$ of that setting. With the input function $\Psi$ and the output functions $\psi_s$, the pointwise assumption of \cite[(5.2)]{CEM} holds with equality. Adding a common constant to the $\psi_s$, and the corresponding constant to $\Psi$, makes these functions nonnegative without changing the conclusion. Hence \cite[Proposition~5.2]{CEM} gives
\begin{equation}\label{eq:cem-bounded}
 \int_{\R^n}\exp\left(-\int_\Omega\psi_s(B_s z)\,d\theta(s)\right)d\gamma_n(z)
 \le
 \exp\int_\Omega\log\left(\int_{\R^k}e^{-\psi_s}\,d\gamma_k\right)d\theta(s).
\end{equation}

The same inequality holds for bounded jointly measurable $\psi_s$. Indeed, apply the preceding argument to $\mathcal Q_t\psi_s$, where $\mathcal Q_t$ is the Ornstein--Uhlenbeck semigroup defined in Section~4. These functions are uniformly bounded and continuous in the spatial variable. As $t\downarrow0$,
\[
 \int_\Omega\int_{\R^k}|\mathcal Q_t\psi_s-\psi_s|\,d\gamma_k\,d\theta(s)\longrightarrow0
\]
by $L^1$ continuity of the semigroup and dominated convergence. Since $B_sZ$ has law $\gamma_k$ for $Z\sim\gamma_n$, the corresponding parameter integrals converge in $L^1(\gamma_n)$. A subsequence converges almost everywhere, and bounded convergence applies to the left side of \eqref{eq:cem-bounded}. The right side converges by the same $L^1$ convergence and the common positive upper and lower bounds of the inner integrals. This proves \eqref{eq:cem-bounded} for bounded measurable families.

Taking $\psi_s=-\log h_s$ gives the bounded positive forward inequality. For a general normalized density family, apply this inequality to $h_s=\min\{g_s,N\}+\varepsilon$ and use the two ordered limits $N\to\infty$ and $\varepsilon\downarrow0$ in the proof of Theorem~\ref{thm:forward} below. Those limits also allow $g_s=0$. Thus the full forward inequality follows from \cite{CEM} together with this approximation.
\end{proof}

\begin{proposition}[Finite parameter spaces and \cite{CDP}]\label{prop:finite-cdp}
Suppose that $\mu$ is supported on finitely many atoms $s_i$ of masses $\mu_i>0$, and put $p_i=p(s_i)$. Then the inequality assertions of Theorems~\ref{thm:forward} and~\ref{thm:reverse} follow from \cite[Theorem~1]{CDP}, including mixed positive and negative exponents and zeros of positive-exponent factors.
\end{proposition}
\begin{proof}
Let $T=(B_iB_j^*)_{i,j}$ and $q_i=p_i/\mu_i$. For a field $h=(h_i)$, set $u_i=\mu_i h_i$. Then
\[
 \langle(C-P)h,h\rangle_{\Hil}
 =\left\langle\left(T-\operatorname{diag}(q_i I_k)\right)u,u\right\rangle.
\]
Thus the two operator comparisons are equivalent to the corresponding matrix comparisons in \cite{CDP}. Apply \cite[Theorem~1]{CDP} to $B_iZ$ with exponents $q_i$ and functions $g_i^{\mu_i/p_i}$. Each $q_i$-moment equals $\int g_i\,d\gamma_k=1$. For $p_i<0$, $g_i$ is positive and finite almost everywhere. There are only finitely many such indices, so their contributions are finite and positive on a common Gaussian full-measure set. On this set, the ordinary finite product equals the product defined by \eqref{eq:forward} or \eqref{eq:reverse}, with value zero whenever a positive-exponent factor vanishes. This proves the asserted inequalities. The compatible-density equality criterion is an additional statement and is not obtained from this application of \cite{CDP}.
\end{proof}

\begin{remark}[A limitation of a direct mixed-sign reduction]\label{rem:cem-mixed}
The reverse condition does not in general allow one to use only the positive-exponent projections in the reconstruction condition $(\mathrm{ii})_1$ of \cite[Setting~2]{CEM}. For example, let $n=2$, $k=1$, and take two unit-mass atoms with
\[
 B_1(x_1,x_2)=x_1,\quad p_1=\tfrac34,\qquad
 B_2(x_1,x_2)=\tfrac{x_1+x_2}{\sqrt2},\quad p_2=-1.
\]
Then
\[
 C-P=\begin{pmatrix}1/4&1/\sqrt2\\[1mm]1/\sqrt2&2\end{pmatrix}\ge0,
 \qquad
 \frac{1}{p_1}B_1^*B_1=\begin{pmatrix}4/3&0\\0&0\end{pmatrix}.
\]
In particular, $B_1e_2=0$, so no linear map applied to the positive projection alone can reconstruct $e_2$. This rules out that particular direct identification, but does not rule out a different reduction of the general mixed-sign inequality to an existing theorem.
\end{remark}
}

{\begin{proposition}[Reduction of the reverse inequality to \cite{CDP}]\label{prop:nonatomic-cdp}
Under the hypotheses of Theorem~\ref{thm:reverse}, its inequality follows from \cite[Theorem~1]{CDP}, finite partition approximations, and the ordered truncation in Section~5. The parameter space may have a nonatomic part, and the positive-exponent factors may vanish.
\end{proposition}

\begin{proof}
We first consider a jointly measurable family $f_s=e^{\phi_s}$ with $|\phi_s|\le M$. We prove
\begin{equation}\label{eq:cdp-reverse-bounded}
 \int_{\R^n}\exp\!\left(\int_\Omega\phi_s(B_s z)\,d\mu(s)\right)d\gamma_n(z)
 \ge
 \exp\!\left(\int_\Omega\frac{1}{p(s)}\log\int_{\R^k}e^{p(s)\phi_s}\,d\gamma_k\,d\mu(s)\right).
\end{equation}
Assume initially that all $\phi_s$ have a common Lipschitz constant. By Lemma~\ref{lem:atoms}, the positive-exponent parameters form finitely many atoms. Include each of them as a singleton in every partition below. On $\{p<0\}$, choose an increasing sequence of finite measurable partitions $\mathcal P_j$ whose generated sigma-algebra is the completed parameter sigma-algebra modulo null sets, and such that the oscillation of $B_s$ in operator norm on every cell is at most $\delta_j$, where $\delta_j\downarrow0$ and $\delta_j<1/2$. Such partitions exist because the parameter space is standard Borel and the range of $B_s$ is bounded in a finite-dimensional space. Null cells are discarded.

{For completeness, choose a countable generator of the Borel sigma-algebra on $\{p<0\}$. At stage $j$, take the common refinement of the first $j$ generating sets and a finite partition obtained by pulling back a finite cover of the range of $B_s$ by sets of diameter at most $\delta_j$. Also refine by all partitions from the preceding stages. The resulting finite partitions are increasing, generate the Borel sigma-algebra, and have the required oscillation. Passing to the measure completion changes their generated sigma-algebra only by null sets. We choose representatives of $p$, $B_s$, and $\phi$ on a common full-measure set before forming the cell averages.}

For a cell $E$ of positive mass $\mu_E$, set
\[
 \overline B_E=\frac1{\mu_E}\int_E B_s\,d\mu(s),\quad
 \overline p_E=\frac1{\mu_E}\int_E p(s)\,d\mu(s),\quad
 \overline\phi_E(x)=\frac1{\mu_E}\int_E\phi_s(x)\,d\mu(s).
\]
For a positive atom, these are its original values. Put $X_E=\overline B_E Z$ and let $T_j$ be the block covariance matrix of the finite vector $(X_E)_{E\in\mathcal P_j}$. For any vectors $u_E\in\R^k$, use the field $h_s=u_E/\mu_E$ on $E$ in $C\ge P$. This gives
\begin{equation}\label{eq:partition-comparison}
 T_j\ge\operatorname{diag}\left(\frac{\overline p_E}{\mu_E}I_k\right)_{E\in\mathcal P_j}.
\end{equation}
If $E\subset\{p<0\}$, choose $s_E\in E$. Since $B_{s_E}B_{s_E}^*=I_k$ and $\|\overline B_E-B_{s_E}\|\le\delta_j$,
\[
 (1-\delta_j)^2I_k\le\overline B_E\overline B_E^*\le I_k.
\]
Set $a_j=(1-\delta_j)^2$ and define nonzero exponents
\[
 q_E=\begin{cases}
 \overline p_E/\mu_E,&E\text{ a positive atom},\\
 \overline p_E/(a_j\mu_E),&E\subset\{p<0\}.
 \end{cases}
\]
Because $\overline p_E<0$ on a negative cell, the preceding covariance bound gives
\[
 q_E\,\operatorname{Cov}(X_E)\le
 \frac{\overline p_E}{\mu_E}I_k.
\]
The same identity holds on every positive atom. Consequently \eqref{eq:partition-comparison} implies the reverse matrix condition in \cite[Theorem~1]{CDP}. Apply that theorem to $X_E$ with exponent $q_E$ and function $x\mapsto e^{\mu_E\overline\phi_E(x)}$. We obtain
\begin{equation}\label{eq:partition-cdp}
 \mathbb E\exp\!\left(\sum_E\mu_E\overline\phi_E(X_E)\right)
 \ge\exp\!\left(\sum_E\frac1{q_E}\log\mathbb E e^{\mu_Eq_E\overline\phi_E(X_E)}\right).
\end{equation}

{In particular, the sign in the negative-cell comparison is essential: if $S_E=\operatorname{Cov}(X_E)\ge a_jI_k$ and $\overline p_E<0$, then
\[
 q_ES_E-\frac{\overline p_E}{\mu_E}I_k
 =\frac{\overline p_E}{a_j\mu_E}(S_E-a_jI_k)\le0.
\]
}

We pass to the limit. Write $\overline B_j(s)$, $\overline p_j(s)$ and $\overline\phi_j(s,x)$ for the corresponding cell averages. The martingale convergence theorem gives
\[
 \overline p_j\to p\quad\text{in }L^1(\mu),\qquad
 \overline\phi_j\to\phi\quad\text{in }L^1(\mu\otimes\gamma_k).
\]
{For the second convergence, condition $\phi(s,x)$ on the product sigma-algebra $\sigma(\mathcal P_j)\otimes\mathcal B(\R^k)$. Its conditional expectation is $\overline\phi_j(s,x)$ almost everywhere. The increasing limit of these sigma-algebras is the product of the completed parameter sigma-algebra with $\mathcal B(\R^k)$, modulo product-null sets. The $L^1$ martingale convergence theorem therefore applies.}
Moreover, $\|\overline B_j(s)-B_s\|\le\delta_j$. If $L$ is the common Lipschitz constant, then
\begin{align*}
 &\mathbb E\left|\sum_E\mu_E\overline\phi_E(\overline B_E Z)
       -\int_\Omega\phi_s(B_sZ)\,d\mu(s)\right|\\
 &\quad\le
 \int_\Omega\int_{\R^k}|\overline\phi_j(s,x)-\phi_s(x)|\,d\gamma_k(x)d\mu(s)
 +L\delta_j\mu(\Omega)\mathbb E|Z|\longrightarrow0.
\end{align*}
Here $B_sZ$ has law $\gamma_k$. Since the exponents are uniformly bounded, the left side of \eqref{eq:partition-cdp} converges to the left side of \eqref{eq:cdp-reverse-bounded}.

For the right side, the covariance matrices $\overline B_E\overline B_E^*$ converge uniformly to $I_k$, so their Gaussian laws converge uniformly in total variation to $\gamma_k$. Indeed, when the eigenvalues of a covariance matrix lie in $[1/4,1]$, its density with respect to Lebesgue measure is bounded by a fixed multiple of $e^{-|x|^2/2}$. The densities converge uniformly on compact sets as the matrices converge to $I_k$, and dominated convergence gives the asserted uniform total variation convergence. On negative cells $\mu_Eq_E=\overline p_E/a_j$ and $1/q_E=\mu_Ea_j/\overline p_E$. On positive atoms $a_j$ is replaced by $1$. Since $|\overline\phi_E|\le M$ and $p_*\le|\overline p_E|\le p^*$, the map
\[
 (r,\psi)\longmapsto\frac1r\log\int_{\R^k}e^{r\psi}\,d\gamma_k
\]
is uniformly continuous in $r$ on $[-p^*,-p_*]\cup[p_*,p^*]$ and is Lipschitz in $\psi$ with respect to $L^1(\gamma_k)$, uniformly for $|\psi|\le M$. The same conclusion holds when $r$ is replaced by the nearby values $\overline p_E/a_j$. Therefore the exponent on the right side of \eqref{eq:partition-cdp} converges to that on the right side of \eqref{eq:cdp-reverse-bounded}. This proves \eqref{eq:cdp-reverse-bounded} for the Lipschitz family.

{More explicitly, for $|r|,|r'|\le4p^*$, $|r|,|r'|\ge p_*/2$, with $rr'>0$, and $|\psi|,|\chi|\le M$, the mean value theorem gives a constant $C=C(M,p_*,p^*)$ such that
\[
 \left|\frac1r\log\int e^{r\psi}\,d\gamma_k
 -\frac1{r'}\log\int e^{r'\chi}\,d\gamma_k\right|
 \le C\left(|r-r'|+\int|\psi-\chi|\,d\gamma_k\right).
\]
The lower bound $\int e^{r\psi}\,d\gamma_k\ge e^{-4p^*M}$ justifies the logarithmic estimate. For sufficiently large $j$ the corrected exponents $\overline p_j/a_j$ lie in these ranges, and
\[
 \int_{\{p<0\}}\left|\frac{\overline p_j}{a_j}-p\right|d\mu
 \le a_j^{-1}\|\overline p_j-p\|_{L^1(\mu)}
   +(a_j^{-1}-1)p^*\mu(\Omega)\longrightarrow0.
\]
Together with the uniform total variation convergence, these estimates justify summing the cell errors even when the number of cells increases.}

For a general bounded measurable $\phi_s$, replace it by $\mathcal Q_t\phi_s$. For each $t>0$ the smoothed family has a common Lipschitz bound, and $\mathcal Q_t\phi_s\to\phi_s$ in $L^1(\mu\otimes\gamma_k)$ as $t\downarrow0$. The coisometry condition, the common bound $M$, and dominated convergence pass both sides of \eqref{eq:cdp-reverse-bounded} to the limit.

Finally apply \eqref{eq:cdp-reverse-bounded} to the bounded positive truncations in \eqref{eq:rev-bounded-approx}. Section~5 proves the successive limits $M\to\infty$, $\varepsilon\downarrow0$, and $N\to\infty$, in this order. Their limit is precisely the extended product $F_-$, even when a positive-exponent factor is zero. Thus the reverse inequality for the stated general class follows. The argument does not imply that the finite approximations are equality cases.
\end{proof}
}

\section{Bounded functions and the Gaussian semigroup}

For $t > 0$, let $\mathcal{Q}_t$ denote the Ornstein--Uhlenbeck semigroup on $\R^d$, defined by
\[
\mathcal{Q}_t h(x) = \int_{\R^d} h\bigl(e^{-t}x + \sqrt{1 - e^{-2t}}\,y\bigr)\,d\gamma_d(y).
\]
(We use calligraphic $\mathcal{Q}_t$ for the semigroup to distinguish it from the orthogonal projections $Q_s = B_s^* B_s$.) With respect to the Gaussian measure $\gamma_d$, $\mathcal{Q}_t$ has the Mehler kernel representation
\[
\mathcal{Q}_t h(x) = \int_{\R^d} h(y) M_t(x,y)\,d\gamma_d(y),
\]
where, setting $a = e^{-t} \in (0,1)$,
\[
M_t(x,y) = (1 - a^2)^{-d/2} \exp\left(\frac{-a^2|y|^2 + 2a\langle x,y\rangle - a^2|x|^2}{2(1 - a^2)}\right).
\]
Completing the square in the exponent gives
\[
\frac{-a^2|y|^2 + 2a\langle x,y\rangle - a^2|x|^2}{2(1 - a^2)} = -\frac{a^2}{2(1-a^2)}\left|y - \frac{x}{a}\right|^2 + \frac{|x|^2}{2}.
\]
Consequently, for any $0 < t_0 \le t \le t_1 < \infty$ and $R > 0$, the kernel $M_t(x,y)$ and its first two $x$-derivatives are uniformly bounded for $|x| \le R$ and $y \in \R^d$. Moreover, for any $R_0 > 0$, $M_t(x,y)$ is bounded below by a positive constant depending only on $t_0, t_1, R, R_0$ whenever $|x| \le R$ and $|y| \le R_0$. In particular, if $h \ge 0$ satisfies $\int_{\R^d} h\,d\gamma_d = 1$ and $\int_{|y| \le R_0} h\,d\gamma_d \ge \eta > 0$, then $\mathcal{Q}_t h(x)$ is bounded above and away from zero uniformly on $[t_0, t_1] \times \{|x| \le R\}$, and the same uniform bounds hold for $|\nabla_x \log \mathcal{Q}_t h(x)|$ and $|\nabla_x^2 \log \mathcal{Q}_t h(x)|$. Furthermore, if $h_j \to h$ in $L^1(\gamma_d)$, then $\mathcal{Q}_t h_j$ and its first two $x$-derivatives converge locally uniformly on $(0,\infty) \times \R^d$ to $\mathcal{Q}_t h$ and its corresponding derivatives.

Suppose now that $(f_s)_{s \in \Omega}$ is a jointly measurable family satisfying $0 < c_0 \le f_s(x) \le C_0 < \infty$ for all $s \in \Omega$ and $x \in \R^k$. For $t > 0$, define
\[
u_s(t,x) = \frac{1}{p(s)} \log \mathcal{Q}_t(f_s^{p(s)})(x), \qquad
U_t(z) = \int_\Omega u_s(t, B_s z)\,d\mu(s), \qquad
\Phi(t) = \log \int_{\R^n} e^{U_t(z)}\,d\gamma_n(z).
\]
Let $L_d = \Delta - x \cdot \nabla$ be the Ornstein--Uhlenbeck generator on $\R^d$. Because $B_s B_s^* = I_k$ and $\partial_t (\mathcal{Q}_t \psi) = L_k(\mathcal{Q}_t \psi)$, the chain rule yields
\[
L_n[u_s(t, B_s z)] = (L_k u_s)(t, B_s z), \qquad
\partial_t u_s = L_k u_s + p(s)|\nabla u_s|^2.
\]
Writing $h_{t,z}(s) = \nabla u_s(t, B_s z)$, we have $\nabla U_t(z) = \int_\Omega B_s^* \nabla u_s(t, B_s z)\,d\mu(s) = A^* h_{t,z}$. Differentiating $\Phi(t)$ for $t > 0$ and integrating by parts with respect to $\gamma_n$ (using $\int_{\R^n} (L_n \varphi)\,e^{U_t}\,d\gamma_n = -\int_{\R^n} \langle \nabla \varphi, \nabla U_t \rangle e^{U_t}\,d\gamma_n$ for $\varphi(z) = U_t(z)$) gives
\begin{align}
\Phi'(t) &= \frac{\int_{\R^n} e^{U_t(z)} \int_\Omega \bigl((L_k u_s)(t, B_s z) + p(s)|\nabla u_s(t, B_s z)|^2\bigr)\,d\mu(s)\,d\gamma_n(z)}{\int_{\R^n} e^{U_t}\,d\gamma_n} \nonumber \\
&= \frac{\int_{\R^n} e^{U_t(z)} \bigl(-\langle \nabla U_t(z), \nabla U_t(z) \rangle + \langle P h_{t,z}, h_{t,z} \rangle_{\Hil}\bigr)\,d\gamma_n(z)}{\int_{\R^n} e^{U_t}\,d\gamma_n} \nonumber \\
&= \frac{\int_{\R^n} e^{U_t(z)} \langle (P - C) h_{t,z}, h_{t,z} \rangle_{\Hil}\,d\gamma_n(z)}{\int_{\R^n} e^{U_t}\,d\gamma_n}. \label{eq:dissipation}
\end{align}
The uniform upper and lower bounds on $f_s$ ensure that $u_s$ and its spatial derivatives up to second order are uniformly bounded on compact subsets of $(0,\infty) \times \R^k$, which justifies differentiation under the integral signs and Gaussian integration by parts.

As $t \downarrow 0$ and $t \to \infty$, bounded convergence yields
\[
\Phi(0) = \log \int_{\R^n} \exp\left(\int_\Omega \log f_s(B_s z)\,d\mu(s)\right) d\gamma_n(z), \qquad
\Phi(\infty) = \int_\Omega \frac{1}{p(s)} \log \left(\int_{\R^k} f_s^{p(s)}\,d\gamma_k\right) d\mu(s).
\]
Consequently, $\Phi'(\cdot) \ge 0$ when $P \ge C$ and $\Phi'(\cdot) \le 0$ when $C \ge P$, establishing both inequalities in the bounded positive class. Moreover, equality holds in either case if and only if $h_{t,z} \in \Ker(P - C)$ for almost every $(t,z) \in (0,\infty) \times \R^n$.

Conversely, the operator conditions $P \ge C$ and $C \ge P$ are necessary for their respective inequalities to hold across all bounded positive families. For $v\in\Hil$ and $R>0$, let
\[
 h_{s,R}(x)=\max\{-R,\min\{R,\langle v(s),x\rangle\}\},
 \qquad f_{s,\varepsilon,R}=e^{\varepsilon h_{s,R}}.
\]
For fixed $R$, the functions $h_{s,R}$ are bounded and odd. Expanding the two endpoints at $\varepsilon=0$ therefore gives
\begin{align*}
 &\lim_{\varepsilon\to0}\frac{2}{\varepsilon^2}
    \bigl(\Phi_{\varepsilon,R}(\infty)-\Phi_{\varepsilon,R}(0)\bigr)\\
 &\quad=\int_\Omega p(s)\int_{\R^k}h_{s,R}(x)^2\,d\gamma_k(x)\,d\mu(s)
 -\int_{\R^n}\left(\int_\Omega h_{s,R}(B_sz)\,d\mu(s)\right)^2d\gamma_n(z).
\end{align*}
The boundedness of $h_{s,R}$ and $p$ justifies this expansion for each fixed $R$. As $R\to\infty$, dominated convergence gives $h_{s,R}(B_s z)\to\langle v(s),B_sz\rangle$ in $L^2(\mu\otimes\gamma_n)$. Consequently the right side converges to
\[
 \int_\Omega p(s)|v(s)|^2\,d\mu(s)
 -\left|\int_\Omega B_s^*v(s)\,d\mu(s)\right|^2
 =\langle(P-C)v,v\rangle_{\Hil}.
\]
The sign of the endpoint difference for every $R$ and $\varepsilon$ therefore forces $P\ge C$ in the forward case and $C\ge P$ in the reverse case. No uniform remainder in $R$ is needed.

\section{Passage to general functions}

We now extend the bounded inequalities to general functions satisfying \eqref{eq:masses}. By scaling $f_s$ by $m_s^{-1/p(s)}$, it suffices to prove the normalized statements in Theorems~\ref{thm:forward} and~\ref{thm:reverse}.

\begin{proof}[Proof of the inequality in Theorem~\ref{thm:forward}]
Assume $p > 0$ almost everywhere and $P \ge C$. For $N \ge 1$ and $\varepsilon > 0$, define $g_{s,N,\varepsilon} = \min\{g_s, N\} + \varepsilon$, which satisfies $\varepsilon \le g_{s,N,\varepsilon} \le N + \varepsilon$. Applying the bounded forward inequality from Section~4 to $f_{s,N,\varepsilon} = g_{s,N,\varepsilon}^{1/p(s)}$ gives
\begin{align*}
	&\int_{\R^n} \exp\left(\int_\Omega \frac{1}{p(s)} \log\bigl(\min\{g_s(B_s z), N\} + \varepsilon\bigr)\,d\mu(s)\right) d\gamma_n(z)
\\&\le \exp\left(\int_\Omega \frac{1}{p(s)} \log\left(\int_{\R^k} (\min\{g_s, N\} + \varepsilon)\,d\gamma_k\right) d\mu(s)\right).
\end{align*}
For fixed $\varepsilon > 0$, letting $N \to \infty$ and applying the monotone convergence theorem (since $\log(\min\{g_s, N\} + \varepsilon) \ge \log \varepsilon$) yields
\[
\int_{\R^n} \exp\left(\int_\Omega \frac{1}{p(s)} \log(g_s(B_s z) + \varepsilon)\,d\mu(s)\right) d\gamma_n(z)
\le \exp\left(\int_\Omega \frac{1}{p(s)} \log(1 + \varepsilon)\,d\mu(s)\right).
\]
Now let $\varepsilon \downarrow 0$. The right-hand side converges to $1$. On the left-hand side, the integrand decreases monotonically to $F_+(z)$ as $\varepsilon \downarrow 0$, and for $\varepsilon \le 1$ it is bounded above by the $\gamma_n$-integrable function corresponding to $\varepsilon = 1$. Dominated convergence gives $\int_{\R^n} F_+\,d\gamma_n \le 1$.
\end{proof}

\begin{proof}[Proof of the inequality in Theorem~\ref{thm:reverse}]
Assume $C \ge P$. Fix an integer $N \ge 1$. For any $M > N^{-1}$ and $\varepsilon > 0$, define
\[
g_{i,N,\varepsilon} = \min\{g_i, N\} + \varepsilon \quad (1 \le i \le m), \qquad
g_{s,N,M} = \min\{\max\{g_s, N^{-1}\}, M\} \quad \text{on } \{p < 0\}.
\]
Since these functions are bounded above and away from zero uniformly in $s$, the bounded reverse inequality from Section~4 applies to $f_s = g_{s,\cdot}^{1/p(s)}$ and yields
\begin{align}
\int_{\R^n} \prod_{i=1}^m g_{i,N,\varepsilon}(B_i z)^{\alpha_i}
&\exp\left(-\int_{\{p < 0\}} \log g_{s,N,M}(B_s z)\,d\nu(s)\right) d\gamma_n(z) \nonumber \\
&\ge \prod_{i=1}^m \left(\int_{\R^k} g_{i,N,\varepsilon}\,d\gamma_k\right)^{\alpha_i}
\exp\left(-\int_{\{p < 0\}} \log\left(\int_{\R^k} g_{s,N,M}\,d\gamma_k\right) d\nu(s)\right). \label{eq:rev-bounded-approx}
\end{align}
Fix $N$ and $\varepsilon \in (0,1)$. On the left-hand side of \eqref{eq:rev-bounded-approx}, the integrand is bounded above by $(N+1)^{\sum_i \alpha_i} N^{\nu(\{p < 0\})}$ and converges pointwise as $M \to \infty$ to the corresponding expression with $\max\{g_s(B_s z), N^{-1}\}$ in place of $g_{s,N,M}(B_s z)$. {To justify the parameter integral in this limit, note that
\[
 -\log N\le\log g_{s,N,M}(B_s z)
 \le\log\bigl(1+g_s(B_s z)\bigr),
\]
and the positive upper bound is integrable on $\{p<0\}\times\R^n$ with respect to $\nu\otimes\gamma_n$. Thus dominated convergence applies for almost every $z$. For the right-hand side, the fiber masses satisfy
\[
 N^{-1}\le\int g_{s,N,M}\,d\gamma_k\le 1+N^{-1},
\]
uniformly in $s$ and $M>N^{-1}$. This bounds their logarithms uniformly and permits a second application of dominated convergence. For fixed $N$ the remaining limit $\varepsilon\downarrow0$ is controlled by the same bounded product.} By dominated convergence on both sides, letting $M \to \infty$ and then $\varepsilon \downarrow 0$ in \eqref{eq:rev-bounded-approx} gives
\[
\int_{\R^n} F_N\,d\gamma_n \ge
\prod_{i=1}^m \left(\int_{\R^k} \min\{g_i, N\}\,d\gamma_k\right)^{\alpha_i}
\exp\left(-\int_{\{p < 0\}} \log\left(\int_{\R^k} \max\{g_s, N^{-1}\}\,d\gamma_k\right) d\nu(s)\right).
\]
As $N \to \infty$, we have $\int_{\R^k} \min\{g_i, N\}\,d\gamma_k \to 1$ for each $i=1,\dots,m$, while
\[
1 \le \int_{\R^k} \max\{g_s, N^{-1}\}\,d\gamma_k \le \int_{\R^k} (g_s + N^{-1})\,d\gamma_k = 1 + N^{-1}.
\]
Thus the right-hand side converges to $1$. Since $F_N \uparrow F_-$ pointwise as $N \to \infty$, the monotone convergence theorem gives $\int_{\R^n} F_-\,d\gamma_n \ge 1$.
\end{proof}

\begin{remark}[Sensitivity to the order of limits]\label{rem:order}
The definition of $F_-$ in \eqref{eq:reverse} depends crucially on taking $M \to \infty$ before $\varepsilon \downarrow 0$. Consider $n = k = 1$ and $B_s = 1$ for all $s$, where $\Omega = \{0\} \cup (0,1)$ with $\mu(\{0\}) = 1$, $p(0) = 1/2$, and $\mu$ equal to Lebesgue measure with $p(s) = -1$ on $(0,1)$. Let
\[
g_0(x) = 2\mathbf{1}_{\{x > 0\}}, \qquad g_s(x) = \exp\left(\frac{x}{s} - \frac{1}{2s^2}\right) \quad (s \in (0,1)).
\]
Every $g_s$ is a probability density on $(\R, \gamma_1)$, and $C \ge P$ holds since $\langle (C-P)h, h\rangle_{\Hil} = (h_0 + \int_0^1 h_s\,ds)^2 - \frac{1}{2}h_0^2 + \int_0^1 h_s^2\,ds \ge 0$ by the Cauchy--Schwarz inequality. Under \eqref{eq:reverse}, we have $F_-(x) = 0$ for all $x \le 0$. However, if one instead replaces $g_0$ by $g_0 + \varepsilon$ for a fixed $\varepsilon > 0$ and removes the lower truncation $N^{-1}$ on $(0,1)$ first, then for $x < 0$ the integral $-\int_0^1 \log g_s(x)\,ds = \int_0^1 (-x/s + 1/(2s^2))\,ds = +\infty$, making the product $+\infty$ for all $x < 0$.
\end{remark}

\begin{example}[A non-Gaussian extremizer with zeros]\label{ex:zeros}
Let $n = 2$ and $k = 1$. Take two atoms of unit mass with
\[
B_1(x_1, x_2) = x_1, \quad p_1 = 1, \qquad B_2(x_1, x_2) = x_2, \quad p_2 = \tfrac{1}{2},
\]
and let $[0,1]$ carry Lebesgue measure with $B_s = B_2$ and $p(s) = -1$. On the first coordinate, the quadratic form of $C - P$ vanishes identically. On the second coordinate, it equals
\[
\left(h_0 + \int_0^1 h_s\,ds\right)^2 - \frac{1}{2} h_0^2 + \int_0^1 h_s^2\,ds \ge 0,
\]
so $C \ge P$. Set $g_1(x) = 2\mathbf{1}_{\{x > 0\}}$ and let $g_2$ be any strictly positive measurable function on $\R$ with $\int_\R g_2\,d\gamma_1 = 1$ and $g_2 < \infty$ almost everywhere. Defining $f_1 = g_1$, $f_2 = g_2^2$, and $f_s = g_2^{-1}$ for $s \in [0,1]$ yields $m_s = 1$ for all $s$ and $F_-(x_1, x_2) = g_1(x_1) g_2(x_2)$. Hence $\int_{\R^2} F_-\,d\gamma_2 = 1$, even though $F_-$ vanishes on the half-space $\{x_1 \le 0\}$.
\end{example}

\section{Equality and regularization}\label{sec:recovery}

{Proposition~\ref{prop:nonatomic-cdp} gives an alternative derivation of the reverse inequality, but it does not transfer equality from the limiting family to its finite approximations. In the argument below, equality is instead preserved by an exact semigroup comparison. The vanishing endpoint defects of bounded approximations are then used in the dissipation identity, and a finite-dimensional closure argument recovers the kernel condition. The subsequent characteristic-function argument returns from positive time to the original densities.}

In this section we prove the necessity of a compatible density when equality holds in Theorem~\ref{thm:forward} or Theorem~\ref{thm:reverse}. Throughout, write $F$ for either $F_+$ or $F_-$, and let $F_t$ denote the corresponding product constructed from the smoothed densities $(\mathcal{Q}_t g_s)_{s \in \Omega}$ for $t > 0$. Since $\mathcal{Q}_t$ is positivity improving, $\mathcal{Q}_t g_s(x) \in (0,\infty)$ for all $t > 0$, $s \in \Omega$, and $x \in \R^k$.

\begin{lemma}[Semigroup comparison]\label{lem:comparison}
For every $t > 0$, we have $\gamma_n$-almost everywhere:
\begin{equation}\label{eq:comparison}
\mathcal{Q}_t F \le F_t \quad \text{in the forward case}, \qquad
F_t \le \mathcal{Q}_t F \quad \text{in the reverse case}.
\end{equation}
\end{lemma}

\begin{proof}
Fix $t > 0$ and $z \in \R^n$, and set $a = e^{-t}$ and $\sigma = \sqrt{1 - e^{-2t}}$. Consider the conditional Gaussian random vector $Z_{t,z} = a z + \sigma \tilde{Z}$ where $\tilde{Z} \sim \gamma_n$, so that $\mathcal{Q}_t \psi(z) = \E[\psi(Z_{t,z})]$ and $\mathcal{Q}_t h(B_s z) = \E[h(B_s Z_{t,z})]$.

In the forward case, apply the bounded forward inequality to the conditional law of $Z_{t,z}$ (equivalently, apply the forward inequality to the functions $x \mapsto \min\{g_s(a B_s z + \sigma x), N\} + \varepsilon$). Letting $N \to \infty$ by monotone convergence gives
\[
\mathcal{Q}_t \left[\exp\left(\int_\Omega \frac{1}{p(s)} \log(g_s(B_s \cdot) + \varepsilon)\,d\mu(s)\right)\right](z)
\le \exp\left(\int_\Omega \frac{1}{p(s)} \log(\mathcal{Q}_t g_s(B_s z) + \varepsilon)\,d\mu(s)\right).
\]
Now let $\varepsilon \downarrow 0$. For $\varepsilon \in (0,1)$, the right-hand side is bounded above by the product formed from $\mathcal{Q}_t g_s + 1$, whose $\gamma_n$-integral is finite by Theorem~\ref{thm:forward} applied to $\mathcal{Q}_t g_s + 1$. Thus for $\gamma_n$-almost every $z$, the integrand inside the exponential on the right-hand side has a $\mu$-integrable positive part, so monotone convergence yields $\mathcal{Q}_t F_+(z) \le F_t(z)$.

In the reverse case, fix an integer $N \ge 1$ and put $g_{i,N} = \min\{g_i, N\}$ for $1 \le i \le m$ and $g_{s,N} = \max\{g_s, N^{-1}\}$ on $\{p < 0\}$. For $M > N^{-1}$ and $\varepsilon > 0$, applying the bounded reverse inequality to the conditional Gaussian vector $Z_{t,z}$ gives
\[
\prod_{i=1}^m \bigl(\mathcal{Q}_t(g_{i,N} + \varepsilon)(B_i z)\bigr)^{\alpha_i}
\exp\left(-\int_{\{p < 0\}} \log \mathcal{Q}_t \min\{g_{s,N}, M\}(B_s z)\,d\nu(s)\right)
\le \mathcal{Q}_t F_{N,M,\varepsilon}(z),
\]
where $F_{N,M,\varepsilon}$ denotes the bounded product in \eqref{eq:rev-bounded-approx}. For fixed $N$ and $\varepsilon \in (0,1)$, we have $\mathcal{Q}_t \min\{g_{s,N}, M\}(B_s z) \ge N^{-1}$ and $\log^+ \mathcal{Q}_t \min\{g_{s,N}, M\}(B_s z) \le \log(1 + \mathcal{Q}_t g_s(B_s z))$, which is $\nu$-integrable for $\gamma_n$-almost every $z$ since $\int_{\R^n} \mathcal{Q}_t g_s(B_s z)\,d\gamma_n(z) = \int_{\R^k} g_s\,d\gamma_k = 1$. Letting $M \to \infty$ and then $\varepsilon \downarrow 0$ via dominated convergence yields
\[
\prod_{i=1}^m \bigl(\mathcal{Q}_t g_{i,N}(B_i z)\bigr)^{\alpha_i}
\exp\left(-\int_{\{p < 0\}} \log \mathcal{Q}_t g_{s,N}(B_s z)\,d\nu(s)\right)
\le \mathcal{Q}_t F_N(z)
\]
for $\gamma_n$-almost every $z$. Finally, let $N \to \infty$. On the left-hand side, $\mathcal{Q}_t g_{i,N}(B_i z) \uparrow \mathcal{Q}_t g_i(B_i z) > 0$, while $\mathcal{Q}_t g_{s,N}(B_s z) \downarrow \mathcal{Q}_t g_s(B_s z)$. Since $\mathcal{Q}_t g_{s,N} \le \mathcal{Q}_t g_s + 1$, the positive parts $\log^+ \mathcal{Q}_t g_{s,N}(B_s z)$ are dominated by $\log(1 + \mathcal{Q}_t g_s(B_s z)) \in L^1(\nu)$ for $\gamma_n$-almost every $z$, while the negative parts $\log^- \mathcal{Q}_t g_{s,N}(B_s z)$ increase monotonically to $\log^- \mathcal{Q}_t g_s(B_s z)$. By monotone convergence, the left-hand side converges to $F_t(z)$, and the right-hand side increases to $\mathcal{Q}_t F_-(z)$.
\end{proof}

Suppose now that equality holds in Theorem~\ref{thm:forward} or Theorem~\ref{thm:reverse}, so that $\int_{\R^n} F\,d\gamma_n = 1$. Since $\mathcal{Q}_t$ preserves $\gamma_n$ and each $\mathcal{Q}_t g_s$ has unit $\gamma_k$-integral, Theorems~\ref{thm:forward} and~\ref{thm:reverse} applied to $(\mathcal{Q}_t g_s)_{s \in \Omega}$ together with Lemma~\ref{lem:comparison} imply that for every $t > 0$,
\begin{equation}\label{eq:exact}
F_t = \mathcal{Q}_t F \quad \gamma_n\text{-almost everywhere}, \qquad \int_{\R^n} F_t\,d\gamma_n = 1.
\end{equation}
In particular, $F_t$ agrees almost everywhere with the strictly positive smooth function $\mathcal{Q}_t F$. The functions involved are jointly measurable in $(t,z)$: this follows from the Mehler formula, the joint measurability of the density family, and the countable truncations in the definition of the product. Fubini's theorem therefore permits us to use the equality simultaneously for almost every $(t,z)$. Notice that this argument does not assume that any of the truncated approximations are extremizers.

\begin{lemma}[Recovery of the kernel identity]\label{lem:kernel-recovery}
If equality holds in Theorem~\ref{thm:forward} or Theorem~\ref{thm:reverse}, then for almost every $(t,z) \in (0,\infty) \times \R^n$,
\begin{equation}\label{eq:field-kernel}
s \longmapsto \frac{1}{p(s)} \nabla \log \mathcal{Q}_t g_s(B_s z) \in \Ker(P - C).
\end{equation}
\end{lemma}

\begin{proof}
We construct a diagonal sequence of bounded positive approximations whose endpoint defects vanish and whose smoothed products converge to $F_t$. We detail the reverse case. The forward case follows identically using $\min\{g_s, N\} + \varepsilon$.

For integers $N \ge 1$, $M > N^{-1}$, and $\varepsilon \in (0,1)$, let $L_{N,M,\varepsilon}$ and $R_{N,M,\varepsilon}$ denote the left- and right-hand sides of \eqref{eq:rev-bounded-approx}, and let $F_{N,M,\varepsilon,t}$ be the smoothed product formed from $\mathcal{Q}_t g_{i,N,\varepsilon}$ and $\mathcal{Q}_t g_{s,N,M}$. As shown in Section~5, $\lim_{N\to\infty} \lim_{\varepsilon\downarrow 0} \lim_{M\to\infty} L_{N,M,\varepsilon} = \int_{\R^n} F_-\,d\gamma_n = 1$ and $\lim_{N\to\infty} \lim_{\varepsilon\downarrow 0} \lim_{M\to\infty} R_{N,M,\varepsilon} = 1$. Moreover, by the same dominated and monotone convergence arguments as in Lemma~\ref{lem:comparison}, for almost every $(t,z) \in (0,\infty) \times \R^n$ we have
\[
\lim_{N\to\infty} \lim_{\varepsilon\downarrow 0} \lim_{M\to\infty} F_{N,M,\varepsilon,t}(z) = F_t(z).
\]
For each integer $j \ge 1$, consider the compact cylinder $E_j = [1/j, j] \times \{z \in \R^n : |z| \le j\}$. Using convergence in $(dt \otimes d\gamma_n)$-measure on $E_j$ successively in $M \to \infty$, $\varepsilon \downarrow 0$, and $N \to \infty$, we may select parameters $(N_j, M_j, \varepsilon_j)$ such that {$N_j\ge j$, $M_j\ge j$, and $0<\varepsilon_j\le j^{-1}$, and such that}:
\begin{enumerate}
\item $|L_{N_j,M_j,\varepsilon_j} - 1| + |R_{N_j,M_j,\varepsilon_j} - 1| < 2^{-j}$;
\item $(dt \otimes d\gamma_n)\bigl(\{(t,z) \in E_j : |F_{N_j,M_j,\varepsilon_j,t}(z) - F_t(z)| \ge 2^{-j}\}\bigr) < 2^{-j}$.
\end{enumerate}
Set $F_{j,t} = F_{N_j,M_j,\varepsilon_j,t}$, $L_j = L_{N_j,M_j,\varepsilon_j}$, $R_j = R_{N_j,M_j,\varepsilon_j}$, and let $h_j(t,z) \in \Hil$ be the corresponding gradient field $s \mapsto p(s)^{-1} \nabla \log \mathcal{Q}_t g_{s,j}(B_s z)$. {For every fixed cylinder $E_\ell$, we have $E_\ell\subset E_j$ for $j\ge\ell$. The second condition and the Borel--Cantelli lemma therefore give $F_{j,t}(z)\to F_t(z)$ for almost every $(t,z)\in E_\ell$. Taking the countable union over $\ell$ proves convergence almost everywhere on $(0,\infty)\times\R^n$.

The same sequence has the required convergence of the individual densities. Indeed, on $\{p<0\}$,
\[
 |g_{s,j}-g_s|\le (g_s-M_j)_++N_j^{-1},
\]
and on each positive atom,
\[
 |g_{i,j}-g_i|\le (g_i-N_j)_++\varepsilon_j.
\]
Since $N_j,M_j\to\infty$ and $\varepsilon_j\to0$, dominated convergence gives $g_{s,j}\to g_s$ in $L^1(\gamma_k)$ for $\mu$-almost every $s$. In fact, integration over the finite measure $\mu$ also gives convergence in $L^1(\mu\otimes\gamma_k)$.}

Let $T \ge 0$ be the self-adjoint operator from Lemma~\ref{lem:gap}. Integrating the dissipation identity \eqref{eq:dissipation} over any compact interval $[a,b] \subset (0,\infty)$ gives
\[
\int_a^b \frac{\int_{\R^n} F_{j,t}(z) \langle T h_j(t,z), h_j(t,z) \rangle_{\Hil}\,d\gamma_n(z)}{\int_{\R^n} F_{j,t}\,d\gamma_n}\,dt
\le |\log L_j - \log R_j| \longrightarrow 0 \quad (j \to \infty).
\]
By monotonicity of $t \mapsto \int_{\R^n} F_{j,t}\,d\gamma_n$, this denominator lies between $L_j$ and $R_j$ and thus stays bounded near $1$. {Hence, on every $E_\ell$, the nonnegative functions
\[
 (t,z)\longmapsto F_{j,t}(z)\langle T h_j(t,z),h_j(t,z)\rangle_{\Hil}
\]
converge to zero in $L^1(dt\otimes\gamma_n)$. A diagonal subsequence over $\ell=1,2,\ldots$ converges to zero almost everywhere on every $E_\ell$. On the same subsequence the already established convergence $F_{j,t}\to F_t$ remains valid. Thus}
\[
F_{j,t}(z) \langle T h_j(t,z), h_j(t,z) \rangle_{\Hil} \longrightarrow 0 \quad \text{for almost every } (t,z) \in (0,\infty) \times \R^n.
\]
Since $F_{j,t}(z) \to F_t(z) = \mathcal{Q}_t F(z) > 0$ almost everywhere, Lemma~\ref{lem:gap} implies
\[
\operatorname{dist}(h_j(t,z), \Ker T) \longrightarrow 0 \quad \text{for almost every } (t,z) \in (0,\infty) \times \R^n.
\]
Fix such a pair $(t,z)$. For every $s \in \Omega$, the densities $g_{s,j}$ converge to $g_s$ in $L^1(\gamma_k)$ as $j \to \infty$, so $\nabla \log \mathcal{Q}_t g_{s,j}(B_s z) \to \nabla \log \mathcal{Q}_t g_s(B_s z)$ for $\mu$-almost every $s$. Since $\mu$ is a finite measure, $h_j(t,z)$ converges in $\mu$-measure to the field in \eqref{eq:field-kernel}. {Let $w_j$ be the orthogonal projection of $h_j(t,z)$ onto $\Ker T$. Then $\|h_j-w_j\|_2\to0$, so $w_j$ has the same limit in $\mu$-measure. The finite-dimensional closure property proved in Section~3 shows that this limit belongs to $\Ker T$. 

For the forward case, use $g_{s,N,\varepsilon}=\min\{g_s,N\}+\varepsilon$. First let $N\to\infty$ and then $\varepsilon\downarrow0$, as in the proof of the forward inequality. Equality makes both endpoint limits equal to one. The conditional comparison gives the corresponding convergence of the smoothed products to $F_t$. On each $E_j$, choose $\varepsilon_j\le j^{-1}$ and then $N_j\ge j$ so that the analogues of the two displayed conditions above hold. Moreover,
\[
 |g_{s,N_j,\varepsilon_j}-g_s|
 \le (g_s-N_j)_+ +\varepsilon_j,
\]
so the densities converge in $L^1(\gamma_k)$ for almost every $s$. The same dissipation and finite-dimensional closure arguments now prove the forward case.}
\end{proof}

\begin{lemma}[A common smooth potential]\label{lem:potential}
Under the hypotheses of Lemma~\ref{lem:kernel-recovery}, set $U_t = \log \mathcal{Q}_t F$ for $t > 0$. Then for every $t > 0$, there exists a set $S_t \subset \Omega$ of full $\mu$-measure such that for all $s \in S_t$ and all $z \in \R^n$,
\begin{equation}\label{eq:common-gradient}
\nabla \log \mathcal{Q}_t g_s(B_s z) = B_s \nabla U_t(z), \qquad \nabla U_t(z) \in V.
\end{equation}
Consequently, $U_t$ splits additively across $Q_s \R^n$ and $(I - Q_s)\R^n$, and the pushforward of $(\mathcal{Q}_t F)\gamma_n$ under $B_s$ equals $(\mathcal{Q}_t g_s)\gamma_k$ for all $s \in S_t$.
\end{lemma}

\begin{proof}
First fix $t_0 > 0$ for which Lemma~\ref{lem:kernel-recovery} holds for $\gamma_n$-almost every $z$. By \eqref{eq:kernel}, for $\gamma_n$-almost every $z$ there exists a unique vector $v(z) \in V$ such that $\nabla \log \mathcal{Q}_{t_0} g_s(B_s z) = B_s v(z)$ for $\mu$-almost every $s$. The vector $v(z)$ may be chosen measurably: on the set where the kernel identity holds, it is equal to $A^*h_z$, where $h_z(s)=p(s)^{-1}\nabla\log\mathcal Q_{t_0}g_s(B_sz)$; set it equal to zero elsewhere. The joint measurability of the integrand and Cauchy--Schwarz justify the parameter integral defining $A^*h_z$. By Fubini's theorem, there is a measurable set $S_{t_0} \subset \Omega$ with $\mu(\Omega \setminus S_{t_0}) = 0$ such that for every $s \in S_{t_0}$, this identity holds for $\gamma_n$-almost every $z$. If $v \in V$ satisfies $B_s v = 0$ for all $s \in S_{t_0}$, then $v = Kv = \int_{S_{t_0}} p(s)^{-1} B_s^* B_s v\,d\mu(s) = 0$. Because $V$ is finite-dimensional, we can choose finitely many points $s_1, \dots, s_d \in S_{t_0}$ such that the linear map $L : V \to (\R^k)^d$ given by $Lv = (B_{s_1} v, \dots, B_{s_d} v)$ is injective. Applying a linear left inverse $L^{-1}$ to the smooth map $z \mapsto (\nabla \log \mathcal{Q}_{t_0} g_{s_j}(B_{s_j} z))_{j=1}^d$ provides a smooth representative of $v : \R^n \to V$. Continuity in $z$ then extends $\nabla \log \mathcal{Q}_{t_0} g_s(B_s z) = B_s v(z)$ to all $z \in \R^n$ for every $s \in S_{t_0}$. If $V=\{0\}$, the same conclusion holds with $v=0$, without choosing any points $s_j$.

Next, we identify $v$ with $\nabla U_{t_0}$. On any compact ball $B \subset \R^n$, the Mehler kernel $M_{t_0}(B_s z, y)$ is uniformly bounded above in $(s,z,y) \in \Omega \times B \times \R^k$. Since $\int_{\R^k} g_s\,d\gamma_k = 1$, the positive parts $\log^+ \mathcal{Q}_{t_0} g_s(B_s z)$ have a uniform upper bound on $\Omega \times B$. On the other hand, by \eqref{eq:exact},
\[
\int_\Omega \frac{1}{p(s)} \log \mathcal{Q}_{t_0} g_s(B_s z)\,d\mu(s) = U_{t_0}(z) \quad \text{for } \gamma_n\text{-a.e. } z,
\]
where $U_{t_0} = \log \mathcal{Q}_{t_0} F$ is smooth and locally bounded. In the forward case ($p > 0$), this forces $\int_B \int_\Omega p(s)^{-1} \log^- \mathcal{Q}_{t_0} g_s(B_s z)\,d\mu(s)\,dz < \infty$. {In the reverse case, the corresponding nonnegative weight on $\{p<0\}$ is $|p(s)|^{-1}$, and the positive-exponent part consists of finitely many smooth terms.} {More explicitly, in the reverse case the contribution of the finitely many positive atoms is bounded on $B$, while the positive logarithmic part on $\{p<0\}$ is uniformly bounded. Since the total signed integral equals the locally bounded function $U_{t_0}$ almost everywhere, the remaining negative logarithmic part has a finite integral over $B\times\{p<0\}$. In the forward case the same conclusion follows without the atomic term. Consequently,
\[
\int_B\int_\Omega\left|\frac{1}{p(s)}\log\mathcal Q_{t_0}g_s(B_s z)\right|\,d\mu(s)\,dz<\infty.
\]
In particular, the parameter integral defines the same distribution as $U_{t_0}$, and Fubini's theorem applies in the calculation below.} Since $\nabla \log \mathcal{Q}_{t_0} g_s(B_s z) = B_s v(z)$ is also locally bounded in $(s,z)$, testing against any $\varphi \in C_c^\infty(\R^n)$ and applying Fubini's theorem gives
\[
-\int_{\R^n} U_{t_0}(z) \nabla \varphi(z)\,dz
= -\int_\Omega \frac{1}{p(s)} \int_{\R^n} \log \mathcal{Q}_{t_0} g_s(B_s z) \nabla \varphi(z)\,dz\,d\mu(s)
= \int_{\R^n} K v(z) \varphi(z)\,dz.
\]
Since $v(z) \in V$, we have $Kv(z) = v(z)$, so $\nabla U_{t_0} = v$ in the sense of distributions and hence everywhere.

We next identify the marginals at this fixed time. For $s\in S_{t_0}$, write $z=B_s^*x+w$ with $w\in\ker B_s$. The identity $\nabla\log\mathcal Q_{t_0}g_s(B_sz)=B_s\nabla U_{t_0}(z)$ implies
\[
 \nabla_x\bigl[U_{t_0}(B_s^*x+w)-\log\mathcal Q_{t_0}g_s(x)\bigr]=0.
\]
Hence there is a positive function $c_s(w)$ such that
\[
 e^{U_{t_0}(B_s^*x+w)}
 =\mathcal Q_{t_0}g_s(x)c_s(w).
\]
Both $e^{U_{t_0}}$ and $\mathcal Q_{t_0}g_s$ have integral one with respect to their standard Gaussian measures. Integrating over $x$ and $w$ therefore gives $\int c_s\,d\gamma_{\ker B_s}=1$, where $\gamma_{\ker B_s}$ is the standard Gaussian probability measure on $\ker B_s$. Thus the $B_s$-marginal of $(\mathcal Q_{t_0}F)\gamma_n$ is $(\mathcal Q_{t_0}g_s)\gamma_k$, and its $Q_s$ and $(I-Q_s)$ components are independent. Moreover, $\nabla U_{t_0}\in V$ makes the density constant in the $V^\perp$ variable, so the $V^\perp$ component is an independent standard Gaussian.

To extend \eqref{eq:common-gradient} from a full-measure set of times $t_0 > 0$ to all $t > 0$, fix any $t > 0$ and pick $t_0 \in (0, t)$ for which \eqref{eq:common-gradient} holds. By the semigroup property $\mathcal{Q}_t = \mathcal{Q}_{t - t_0} \mathcal{Q}_{t_0}$, the density $\mathcal{Q}_t F$ equals $\mathcal{Q}_{t - t_0}(\mathcal{Q}_{t_0} F)$ and $\mathcal{Q}_t g_s = \mathcal{Q}_{t - t_0}(\mathcal{Q}_{t_0} g_s)$. Since $U_{t_0}$ satisfies \eqref{eq:common-gradient}, differentiating gives $Q_s \nabla^2 U_{t_0} = \nabla^2 U_{t_0} Q_s$ and $P_V \nabla^2 U_{t_0} = \nabla^2 U_{t_0} P_V$, so $\mathcal{Q}_{t_0} F$ factors as a product across $Q_s \R^n \oplus (I - Q_s)\R^n$ and across $V \oplus V^\perp$. Convolution with the product Gaussian kernel in $\mathcal{Q}_{t - t_0}$ preserves these product structures, so $U_t = \log \mathcal{Q}_t F$ also splits additively across $Q_s \R^n \oplus (I - Q_s)\R^n$ and across $V \oplus V^\perp$. The marginal identity at time $t_0$, followed by the semigroup property, shows that the $B_s$-marginal at time $t$ is $(\mathcal Q_tg_s)\gamma_k$. The preservation of independence then gives $\nabla\log\mathcal Q_tg_s(B_sz)=B_s\nabla U_t(z)$ for all $z$. The standard Gaussian factor on $V^\perp$ is also preserved, and hence $\nabla U_t(z)\in V$.
\end{proof}

\begin{proof}[Proof of necessity in Theorems~\ref{thm:forward} and~\ref{thm:reverse}]
Suppose equality holds in either theorem, so $$\int_{\R^n} F\,d\gamma_n = 1.$$ Let $Y$ be a random vector with law $F \gamma_n$. Fix $t = 1$ and set $a = e^{-1}$, $\sigma = \sqrt{1 - e^{-2}}$. If $\tilde{Z} \sim \gamma_n$ is independent of $Y$, then $Y_1 = a Y + \sigma \tilde{Z}$ has density $\mathcal{Q}_1 F$ with respect to $\gamma_n$.

For any $\xi \in \R^n$, the characteristic function of $Y_1$ satisfies
\[
\E[e^{i\langle \xi, Y_1 \rangle}] = \E[e^{i\langle a\xi, Y \rangle}] \exp\left(-\frac{\sigma^2 |\xi|^2}{2}\right).
\]
By Lemma~\ref{lem:potential}, for every $s \in S_1$, $B_s Y_1$ has law $(\mathcal{Q}_1 g_s)\gamma_k$. On the other hand, $B_s Y_1 = a B_s Y + \sigma B_s \tilde{Z}$ where $B_s \tilde{Z} \sim \gamma_k$ is independent of $B_s Y$, so for every $\eta \in \R^k$,
\begin{align*}
	\E[e^{i\langle a\eta, B_s Y \rangle}] \exp\left(-\frac{\sigma^2 |\eta|^2}{2}\right)
&= \int_{\R^k} e^{i\langle \eta, x \rangle} (\mathcal{Q}_1 g_s)(x)\,d\gamma_k(x)
\\&= \left(\int_{\R^k} e^{i\langle a\eta, x \rangle} g_s(x)\,d\gamma_k(x)\right) \exp\left(-\frac{\sigma^2 |\eta|^2}{2}\right).
\end{align*}
Cancelling the strictly positive Gaussian factor $\exp(-\sigma^2|\eta|^2/2)$ shows that $B_s Y$ has law $g_s \gamma_k$ for all $s \in S_1$.

Similarly, for any $s \in S_1$ and any $\xi_1 \in Q_s \R^n$, $\xi_2 \in (I - Q_s)\R^n$, the independence of $Q_s Y_1$ and $(I - Q_s)Y_1$ gives
\[
\E[e^{i\langle a(\xi_1 + \xi_2), Y \rangle}] e^{-\sigma^2(|\xi_1|^2 + |\xi_2|^2)/2}
= \E[e^{i\langle a\xi_1, Q_s Y \rangle}] e^{-\sigma^2|\xi_1|^2/2} \,
  \E[e^{i\langle a\xi_2, (I - Q_s)Y \rangle}] e^{-\sigma^2|\xi_2|^2/2}.
\]
Cancelling the nonzero Gaussian factors establishes the independence of $Q_s Y$ and $(I - Q_s)Y$. Finally, since $\nabla U_1(z) \in V$ for all $z$, the density $\mathcal{Q}_1 F$ is constant in the $V^\perp$ variables, meaning that $P_V Y_1$ and $(I - P_V)Y_1$ are independent and $(I - P_V)Y_1$ is standard Gaussian on $V^\perp$. Applying the same characteristic-function cancellation to $\xi = \xi_V + \xi_{V^\perp}$ shows that $P_V Y$ and $(I - P_V)Y$ are independent and that $\E[e^{i\langle a\xi_{V^\perp}, (I - P_V)Y \rangle}] = \exp(-a^2 |\xi_{V^\perp}|^2/2)$, so $(I - P_V)Y$ is standard Gaussian on $V^\perp$. Thus $r = F$ is a compatible probability density.
\end{proof}

\section{Support and sufficiency}\label{sec:support}

\begin{lemma}[Support of a compatible density]\label{lem:support}
In the reverse setting $C \ge P$, let $r$ be a compatible probability density on $(\R^n, \gamma_n)$. Then
\begin{equation}\label{eq:support}
\{z \in \R^n : r(z) = 0\} = \bigcup_{i=1}^m \{z \in \R^n : g_i(B_i z) = 0\}
\quad \text{modulo } \gamma_n\text{-null sets},
\end{equation}
where the union is empty if $m = 0$.
\end{lemma}

\begin{proof}
{Throughout this proof, assertions involving the indices $s_i$ are understood after choosing representatives on a common full-measure set.} Since $B_i Y$ has density $g_i$ with respect to $\gamma_k$ when $Y \sim r\gamma_n$, we have $\int_{\{g_i(B_i \cdot) = 0\}} r\,d\gamma_n = \int_{\{g_i = 0\}} g_i\,d\gamma_k = 0$, which gives the inclusion from right to left in \eqref{eq:support}.

For the reverse inclusion, let $F_1, \dots, F_\ell$ be the nonzero common eigenspaces of the orthogonal projections $Q_1, \dots, Q_m$ (corresponding to eigenvalues in $\{0,1\}$), and set $D = (\bigoplus_{j=1}^\ell F_j)^\perp$. Thus $\R^n = D \oplus F_1 \oplus \dots \oplus F_\ell$ is an orthogonal decomposition, and each $Q_i$ restricts to either $0$ or the identity on each $F_j$. Call a block $F_j$ \emph{visible} if $Q_i|_{F_j} = I$ for at least one $i \in \{1,\dots,m\}$, and \emph{invisible} if $Q_i|_{F_j} = 0$ for all $i = 1,\dots,m$.

Fix $t > 0$ and set $r_t = \mathcal{Q}_t r$. Since compatibility is preserved under Gaussian convolution, {$Q_sY_t$ and $(I - Q_s)Y_t$ are independent} for $\mu$-almost every $s$ (in particular, for every atom $s_1, \dots, s_m$), and $(I - P_V)Y_t$ is standard Gaussian on $V^\perp$ for $Y_t \sim r_t \gamma_n$. Consequently, $\log r_t$ splits additively across $Q_i \R^n \oplus (I - Q_i)\R^n$ and across $V \oplus V^\perp$, which implies that its Hessian commutes with every $Q_i$ and satisfies $\nabla \log r_t(z) \in V$ for all $z \in \R^n$.

We claim that $\log r_t$ admits the block decomposition
\begin{equation}\label{eq:block-decomp}
\log r_t(z) = c + \frac{1}{2}\langle H z_D, z_D \rangle + \langle b, z_D \rangle + \sum_{j=1}^\ell \psi_j(z_j),
\end{equation}
where $z = z_D + \sum_{j=1}^\ell z_j$ is the orthogonal decomposition of $z$. Indeed, for a fixed $z \in \R^n$ and $v \in \R^n$, let $C_v = \partial_v (\nabla^2 \log r_t)(z)$ be the directional derivative of the Hessian in direction $v$. By symmetry of third derivatives and differentiation of $Q_i \nabla^2 \log r_t = (\nabla^2 \log r_t) Q_i$, we have $C_v w = C_w v$ and $Q_i C_v = C_v Q_i$ for all $v, w \in \R^n$ and $1 \le i \le m$. Hence, for any $v, w \in \R^n$,
\[
C_{Q_i v} w = C_w (Q_i v) = Q_i C_w v = Q_i C_v w,
\]
so $C_{Q_i v} = Q_i C_v$. Because $C_{Q_i v}$ also commutes with $Q_j$, it follows that
\[
Q_i Q_j C_v = Q_i C_{Q_j v} = C_{Q_j v} Q_i = Q_j C_v Q_i = Q_j Q_i C_v.
\]
Define $W = \operatorname{span}\{C_v w : v, w \in \R^n\}$. The relation $Q_i C_v = C_{Q_i v}$ implies $Q_i(W) \subset W$, and the previous identity shows that the restrictions $Q_1|_W, \dots, Q_m|_W$ commute. Since they are commuting orthogonal projections on $W$, they share a common orthonormal eigenbasis with eigenvalues in $\{0,1\}$, whence $W \subset \bigoplus_{j=1}^\ell F_j$. As each $C_v$ is symmetric with range in $W$, we have $C_v|_D = 0$ for all $v \in \R^n$. Thus the $D \times D$ block of $\nabla^2 \log r_t(z)$ is constant in $z$. {The commutation relations also show that $\nabla^2\log r_t(z)$ leaves every complete common eigenspace $F_j$ invariant: if $u\in F_j$ and $Q_i u=\varepsilon_i u$, then $Q_i(\nabla^2\log r_t(z)u)=\varepsilon_i\nabla^2\log r_t(z)u$ for every $i$. Symmetry then makes $D$ invariant as well. Thus all mixed Hessian blocks between the spaces $D,F_1,\ldots,F_\ell$ vanish. Integrating along these mutually orthogonal spaces gives a sum of functions of the separate variables. The $D$-function is quadratic because its Hessian is constant. This proves \eqref{eq:block-decomp}.}

If $F_j$ is invisible and $v \in F_j \cap V$, then $Q_i v = 0$ for all $1 \le i \le m$, so
\[
|v|^2 = \langle K v, v \rangle = -\int_{\{p < 0\}} \frac{|Q_s v|^2}{|p(s)|}\,d\mu(s) \le 0,
\]
giving $F_j \cap V = \{0\}$. Since $\nabla \log r_t(z) \in V$, varying $z_j$ in \eqref{eq:block-decomp} gives, for all $u,v\in F_j$,
\[
 \nabla\psi_j(u)-\nabla\psi_j(v)
 =\nabla\log r_t(z_D+u+\textstyle\sum_{h\ne j}z_h)
 -\nabla\log r_t(z_D+v+\textstyle\sum_{h\ne j}z_h).
\]
The difference belongs to $F_j\cap V=\{0\}$. Thus $\nabla\psi_j$ is constant, and $\psi_j$ is affine on every invisible block. The constant gradient need not vanish.

Thus $r_t \gamma_n$ factors as a product measure over $D \oplus F_1 \oplus \dots \oplus F_\ell$, {where its $D$-marginal is Gaussian and each invisible $F_j$-marginal is a Gaussian with identity covariance, possibly with a nonzero mean.} {We now fix $t=1$.} Writing $Y_1 = a Y + \sigma \tilde{Z}$ as in Section~6 with $Y \sim r\gamma_n$, we obtain for $\xi = \xi_D + \sum_{j=1}^\ell \xi_j$:
\[
\E[e^{i\langle a\xi, Y \rangle}] e^{-\sigma^2 |\xi|^2/2}
= \E[e^{i\langle \xi_D, Y_{1,D} \rangle}] \prod_{j=1}^\ell \E[e^{i\langle \xi_j, Y_{1,j} \rangle}].
\]
Setting all but one component of $\xi$ to zero and substituting back allows us to cancel the nonzero Gaussian factors $\exp(-\sigma^2|\xi|^2/2)$, proving that $Y$ itself factors independently across $D \oplus F_1 \oplus \dots \oplus F_\ell$. {For the $D$-component, let $m_{D,1}$ and $\Sigma_{D,1}$ be the mean and covariance of the Gaussian vector $Y_{1,D}$. Division by the Gaussian characteristic function gives
\[
 \mathbb E e^{i\langle\xi_D,Y_D\rangle}
 =\exp\!\left(i\left\langle \frac{m_{D,1}}a,\xi_D\right\rangle
 -\frac12\left\langle\frac{\Sigma_{D,1}-\sigma^2 I_D}{a^2}\xi_D,\xi_D\right\rangle\right).
\]
Set $\Sigma_D=(\Sigma_{D,1}-\sigma^2 I_D)/a^2$.} Boundedness of the characteristic function forces $\Sigma_D \ge 0$, and absolute continuity of the marginal of $r\gamma_n$ on $D$ prevents $\Sigma_D$ from having a zero eigenvalue (which would concentrate $Y_D$ on a proper affine subspace). Hence $Y_D$ is a nondegenerate Gaussian vector with strictly positive density $r_D$ on $D$. {Similarly, an invisible component $Y_j$ deconvolves from a Gaussian with identity covariance to a Gaussian with identity covariance. Its mean may be nonzero, but its density $r_j$ with respect to the standard Gaussian measure is strictly positive.}

Write $r(z) = r_D(z_D) \prod_{j=1}^\ell r_j(z_j)$ for $\gamma_n$-almost every $z$. For each $i \in \{1,\dots,m\}$, since $Q_i$ acts as either $0$ or $I$ on each $F_j$, integrating $r$ over $(I - Q_i)\R^n$ with respect to the standard Gaussian measure yields
\[
g_i(B_i z) = r_{D,i}(Q_i z_D) \prod_{\{j : Q_i|_{F_j} = I\}} r_j(z_j) \qquad \text{for } \gamma_n\text{-a.e. } z,
\]
where $r_{D,i}$ is the strictly positive Gaussian density of $Q_i Y_D$. Because every visible block $F_j$ satisfies $Q_i|_{F_j} = I$ for at least one $i \in \{1,\dots,m\}$, we conclude that
\[
\{r = 0\} = \bigcup_{j \text{ visible}} \{r_j(z_j) = 0\} = \bigcup_{i=1}^m \{g_i(B_i z) = 0\}
\quad \text{modulo } \gamma_n\text{-null sets}.
\]
Finally, if $m = 0$, then $K \le 0$ and $V = \Ker(I - K) = \{0\}$, so compatibility forces $r = 1$ almost everywhere.
\end{proof}

\begin{proof}[Proof of sufficiency in Theorems~\ref{thm:forward} and~\ref{thm:reverse}]
Let $r$ be a compatible density and set $r_t = \mathcal{Q}_t r$ for $t > 0$. By compatibility and Gaussian convolution, for every $t > 0$ and $\mu$-almost every $s$, the pushforward of $r_t \gamma_n$ under $B_s$ is $(\mathcal{Q}_t g_s)\gamma_k$, and $Q_s Y_t$ is independent of $(I - Q_s)Y_t$ for $Y_t \sim r_t \gamma_n$. Moreover, $(I - P_V)Y_t$ is independent of $P_V Y_t$ and standard Gaussian on $V^\perp$.

The exceptional parameter set can be chosen independently of $t$. Indeed, select a full-measure set on which both compatibility statements involving $B_s$ hold at time zero. For each parameter in this set, the Ornstein--Uhlenbeck evolution acts separately on the two orthogonal components, and the image under $B_s$ evolves by the $k$-dimensional semigroup. Hence the asserted independence and marginal identity hold for every $t>0$ on this same parameter set. Since the smoothed densities are positive and smooth, their almost-everywhere product factorizations have versions valid at every spatial point.

Since all marginals $g_s \gamma_k$ arise as pushforwards of a single probability measure $r \gamma_n$ under contractions $B_s$, the family $\{g_s \gamma_k\}_{s \in \Omega}$ is uniformly tight: for every $\eta \in (0,1)$, there exists $R_0 > 0$ such that $\int_{|y| \le R_0} g_s(y)\,d\gamma_k(y) \ge \eta$ {for $\mu$-almost every $s\in\Omega$. We restrict to a common full-measure set of parameters.} By the Mehler kernel bounds in Section~4, on any compact subset of $(0,\infty) \times \R^k$, the functions $\mathcal{Q}_t g_s$ are bounded above and away from zero uniformly in $s \in \Omega$, and their spatial and temporal derivatives up to second order are uniformly bounded in $s \in \Omega$.

Independence of $Q_s Y_t$ and $(I - Q_s)Y_t$ implies $\nabla \log \mathcal{Q}_t g_s(B_s z) = B_s \nabla \log r_t(z)$ for all $z \in \R^n$ and $\mu$-almost every $s$, while the standard Gaussian factor on $V^\perp$ gives $\nabla \log r_t(z) \in V$. Integrating against $p(s)^{-1}\,d\mu(s)$ yields
\[
\nabla \int_\Omega \frac{1}{p(s)} \log \mathcal{Q}_t g_s(B_s z)\,d\mu(s)
= K \nabla \log r_t(z) = \nabla \log r_t(z).
\]
Furthermore, since $\partial_t \log \mathcal{Q}_t g_s = L_k \log \mathcal{Q}_t g_s + |\nabla \log \mathcal{Q}_t g_s|^2$ and $B_s B_s^* = I_k$, differentiating in $t$ gives
\begin{align*}
\partial_t \int_\Omega \frac{1}{p(s)} \log \mathcal{Q}_t g_s(B_s z)\,d\mu(s)
&= \int_\Omega \frac{1}{p(s)} \bigl(L_n[\log \mathcal{Q}_t g_s(B_s z)] + |B_s \nabla \log r_t(z)|^2\bigr)\,d\mu(s) \\
&= L_n \log r_t(z) + \langle K \nabla \log r_t(z), \nabla \log r_t(z) \rangle \\
&= L_n \log r_t(z) + |\nabla \log r_t(z)|^2 = \partial_t \log r_t(z).
\end{align*}
Hence $\log F_t(z) - \log r_t(z)$ is constant in both $t > 0$ and $z \in \R^n$. At $z=0$, the Mehler kernel tends to $1$ uniformly on every compact set of the integration variable and is bounded uniformly for all sufficiently large $t$. The uniform tightness of the marginal laws therefore gives $\sup_s|\mathcal Q_tg_s(0)-1|\to0$. The same kernel argument gives $r_t(0)\to1$. Evaluating the constant at $z=0$ and letting $t\to\infty$ shows that it is zero. Thus $F_t=r_t$ for all $t>0$ and $z\in\R^n$.

Now let $t \downarrow 0$, so that $r_t \to r$ in $L^1(\gamma_n)$ and $\mathcal{Q}_t g_s \to g_s$ in $L^1(\gamma_k)$ for $\mu$-almost every $s$. In the forward case, $\log^+ \mathcal{Q}_t g_s(B_s z)$ is uniformly integrable in $(s,z)$ with respect to $\mu \otimes \gamma_n$ because $\log(1+x) = o(x)$ as $x \to \infty$ and $\int_{\R^n} \mathcal{Q}_t g_s(B_s z)\,d\gamma_n(z) = 1$. {Since $\mathcal Q_t g_s\to g_s$ in $L^1(\gamma_k)$ for almost every $s$, dominated convergence gives $\mathcal Q_t g_s(B_s z)\to g_s(B_s z)$ in $L^1(\mu\otimes\gamma_n)$. The inequality $|\log^+x-\log^+y|\le |x-y|$ gives convergence of the positive logarithmic parts in $L^1(\mu\otimes\gamma_n)$. We may therefore choose a sequence $t_j\downarrow0$ such that $r_{t_j}\to r$ almost everywhere, the positive logarithmic parts converge in $L^1(\mu)$ for almost every $z$, and the densities converge for $(\mu\otimes\gamma_n)$-almost every $(s,z)$. Applying Fatou's lemma to the negative logarithmic parts yields}
\[
r(z) = \lim_{j \to \infty} F_{t_j}(z) \le F_+(z) \quad \text{for } \gamma_n\text{-a.e. } z.
\]
Since $\int_{\R^n} r\,d\gamma_n = 1$ and $\int_{\R^n} F_+\,d\gamma_n \le 1$, it follows that $r = F_+$ almost everywhere and $\int_{\R^n} F_+\,d\gamma_n = 1$.

In the reverse case, for $\mu$-almost every $s$, define the complementary marginal density $k_s$ on $(I - Q_s)\R^n$ and its smoothed version $k_{s,t}$ by
\[
k_s((I - Q_s)z) = \int_{\R^n} r\bigl((I - Q_s)z + Q_s y\bigr)\,d\gamma_n(y), \qquad
k_{s,t}((I - Q_s)z) = \int_{\R^n} r_t\bigl((I - Q_s)z + Q_s y\bigr)\,d\gamma_n(y).
\]
Here we use a nonnegative Borel representative of $r$. These expressions are jointly measurable in their displayed variables: $(s,z,y)\mapsto(I-Q_s)z+Q_sy$ is measurable, and integration over $y$ preserves measurability. The smoothed expression is also jointly measurable in $t>0$. The choice of the representative does not change either expression almost everywhere. Indeed, for each fixed $s$, if $Z$ and $Y$ are independent standard Gaussian vectors, then $(I-Q_s)Z+Q_sY$ is again standard Gaussian. Fubini's theorem therefore transfers any Gaussian null exceptional set to a $(\mu\otimes\gamma_n\otimes\gamma_n)$-null set.
By compatibility, $r(z) = g_s(B_s z) k_s((I - Q_s)z)$ and $r_t(z) = \mathcal{Q}_t g_s(B_s z) k_{s,t}((I - Q_s)z)$ for $\gamma_n$-almost every $z$ and $\mu$-almost every $s$. Since conditional expectation is an $L^1$-contraction,
\[
\int_{\{p < 0\}} \int_{\R^n} |k_{s,t}((I - Q_s)z) - k_s((I - Q_s)z)|\,d\gamma_n(z)\,d\nu(s)
\le \nu(\{p < 0\}) \|r_t - r\|_{L^1(\gamma_n)} \longrightarrow 0
\]
as $t \downarrow 0$, with the same $L^1(\nu \otimes \gamma_n)$ convergence holding for $\mathcal{Q}_t g_s(B_s z) \to g_s(B_s z)$. Choose a sequence $t_j \downarrow 0$ such that for $\gamma_n$-almost every $z$, we have $r_{t_j}(z) \to r(z)$, $\mathcal{Q}_{t_j} g_i(B_i z) \to g_i(B_i z)$ for $1 \le i \le m$, and both $\mathcal{Q}_{t_j} g_s(B_s z) \to g_s(B_s z)$ and $k_{s,t_j}((I - Q_s)z) \to k_s((I - Q_s)z)$ in $L^1(\nu)$.

More precisely, choose $t_j$ so that the sum of the relevant $L^1$ errors is at most $2^{-j}$. Tonelli's theorem then gives summability of the corresponding parameter integrals for almost every $z$, and the usual subsequence argument gives the stated pointwise convergences. The factorizations at time zero and at all times $t_j$ hold outside one common $(\nu\otimes\gamma_n)$-null set, because the time set is countable. The finitely many positive atoms can be included by removing one further Gaussian null set.

Fix such a point $z$ with $r(z) > 0$. The marginal factorization $r(z)=g_i(B_i z)k_i((I-Q_i)z)$ for the finitely many positive atoms implies $0<g_i(B_i z)<\infty$ for every $i$, outside a common Gaussian null set. By Fubini, $0<g_s(B_s z)<\infty$ also holds for $\nu$-almost every negative parameter $s$. Thus the logarithms and the positive factors appearing below have their ordinary meanings. Then $r_{t_j}(z) \to r(z) > 0$, and the identity $\mathcal{Q}_{t_j} g_s(B_s z) k_{s,t_j}((I - Q_s)z) = r_{t_j}(z)$ implies
\[
\log^- \mathcal{Q}_{t_j} g_s(B_s z) \le |\log r_{t_j}(z)| + \log^+ k_{s,t_j}((I - Q_s)z).
\]
Because $\log^+ x \le x$ is $1$-Lipschitz on $[0,\infty)$ and $k_{s,t_j}((I - Q_s)z)$ converges in $L^1(\nu)$, the right-hand side is uniformly $\nu$-integrable in $j$, while $\log^+ \mathcal{Q}_{t_j} g_s(B_s z)$ is also uniformly $\nu$-integrable in $j$ by the $L^1(\nu)$ convergence of $\mathcal{Q}_{t_j} g_s(B_s z)$. Convergence in $L^1(\nu)$ implies convergence in measure. Continuity of the logarithm at the positive limit gives convergence in measure of the logarithms. By Vitali's convergence theorem, $\int_{\{p < 0\}} \log \mathcal{Q}_{t_j} g_s(B_s z)\,d\nu(s) \to \int_{\{p < 0\}} \log g_s(B_s z)\,d\nu(s)$. The limiting logarithm is integrable. Since every positive-atom factor is strictly positive and finite at $z$, the ordered truncation defining $F_-$ agrees there with the exponential of this ordinary logarithmic integral. Consequently $F_-(z) = \lim_{j \to \infty} F_{t_j}(z) = r(z)$. On the other hand, for $\gamma_n$-almost every $z$ with $r(z) = 0$, Lemma~\ref{lem:support} gives $g_i(B_i z) = 0$ for some $i \in \{1,\dots,m\}$, so $F_-(z) = 0 = r(z)$ by \eqref{eq:reverse}. Consequently, $F_- = r$ $\gamma_n$-almost everywhere and $\int_{\R^n} F_-\,d\gamma_n = 1$.
\end{proof}

\begin{proof}[Proof of Corollary~\ref{cor:nonatomic}]
{By Lemma~\ref{lem:atoms}, $p<0$ almost everywhere on the nonatomic part of $\mu$. In particular, if $\mu$ is nonatomic, put $d\nu=|p|^{-1}d\mu$ and
\[
 H(z)=-\int_\Omega\log g_s(B_s z)\,d\nu(s).
\]
The positive part of $\log g_s(B_s z)$ is integrable on $\Omega\times\R^n$, since $\log^+u\le u$ and every $g_s$ has unit mass. Thus $H$ is well defined almost everywhere, $H^-$ is integrable, and $F_-=e^H$. Jensen's inequality, first in the Gaussian variable on each fiber and then in $z$, gives in the extended sense
\[
 \int_\Omega\int_{\R^k}\log g_s\,d\gamma_k\,d\nu(s)\le0,
 \qquad
 \int_{\R^n}F_-\,d\gamma_n
 \ge \exp\!\left(\int_{\R^n}H\,d\gamma_n\right)\ge1.
\]
If equality holds, both displayed inequalities are equalities. In particular, $\int\log g_s\,d\gamma_k=0$ for $\nu$-almost every $s$. Strict concavity of $\log$ and $\int g_s\,d\gamma_k=1$ imply $g_s=1$ almost everywhere for $\nu$-almost every $s$. The converse is immediate. This also gives a direct proof of the purely nonatomic lower bound and its equality case using Jensen's inequality.}
\end{proof}

\section{An explicit structure of the equality cases}\label{sec:explicit}

We now describe compatible densities without using an unknown joint density in the statement. Let $\mathcal L$ be the smallest linear subspace of the symmetric operators on $\R^n$ containing $Q_s$ for $\mu$-almost every $s$. It exists: choose a full-measure set for which the dimension of the span of its projections is minimal. Removing any further null set cannot decrease this dimension, so every other full-measure set has a span containing this minimal span. Choose a full-measure set $S_0$ on which $B_sB_s^*=I_k$ and $Q_s\in\mathcal L$. Its projections span $\mathcal L$. Choose finitely many parameters $s_1',\ldots,s_d'\in S_0$ whose projections span $\mathcal L$, and let $F_1,\ldots,F_\ell$ be the nonzero complete common eigenspaces of $Q_{s_1'},\ldots,Q_{s_d'}$. If $\mathcal L=\{0\}$, take the single common eigenspace $\R^n$. The spaces $F_j$ are pairwise orthogonal. Since every $Q_s$, $s\in S_0$, is a linear combination of the selected projections, its restriction to $F_j$ is a scalar orthogonal projection. Hence it is either zero or the identity. Put
\begin{equation}\label{eq:explicit-spaces}
 A_j=F_j\cap V,\qquad
 A=\bigoplus_{j=1}^\ell A_j,\qquad H=A^\perp,\qquad V_H=V\cap H.
\end{equation}
The zero spaces $A_j$ may be omitted. The subspaces in \eqref{eq:explicit-spaces} depend only on $\mathcal L$ and $V$, apart from the labeling of the complete common eigenspaces. In particular, they are unchanged when $S_0$ is replaced by another full-measure set contained in $\{s:Q_s\in\mathcal L\}$. A full-measure set that also contains exceptional values outside $\mathcal L$ must first be restricted in this way.  For $s\in S_0$, both $A$ and $H$ reduce $Q_s$. Write $Q_{s,H}=Q_s|_H$ and
\[
 J_s=\{j:A_j\ne\{0\},\ Q_s|_{A_j}=I_{A_j}\}.
\]
We write $\gamma_L$ for the standard Gaussian probability measure on a Euclidean subspace $L$, including the probability measure at zero if $L=\{0\}$.

\begin{theorem}[Explicit equality structure]\label{thm:explicit}
Assume the hypotheses of either Theorem~\ref{thm:forward} or Theorem~\ref{thm:reverse}. Equality holds if and only if there are a vector $b\in V_H$, a positive definite symmetric operator $\Sigma$ on $H$, and probability densities $h_j$ with respect to $\gamma_{A_j}$ such that
\begin{equation}\label{eq:explicit-algebra}
 (\Sigma-I_H)H\subset V_H,\qquad
 \Sigma Q_{s,H}=Q_{s,H}\Sigma\quad\text{for }\mu\text{-almost every }s,
\end{equation}
and, for $\mu$-almost every $s$ and $\gamma_n$-almost every $z$,
\begin{equation}\label{eq:explicit-factors}
 g_s(B_sz)=G_s(Q_{s,H}z_H)\prod_{j\in J_s}h_j(z_j),
 \quad z_H=P_Hz,\quad z_j=P_{A_j}z.
\end{equation}
Here $G_s$ is the density with respect to $\gamma_{Q_{s,H}H}$ of the Gaussian distribution on $Q_{s,H}H$ with mean $b_s=Q_{s,H}b$ and covariance $\Sigma_s=Q_{s,H}\Sigma|_{Q_{s,H}H}$. Explicitly,
\[
 G_s(x)=(\det\Sigma_s)^{-1/2}
 \exp\!\left(-\tfrac12\langle\Sigma_s^{-1}(x-b_s),x-b_s\rangle
 +\tfrac12|x|^2\right),\qquad x\in Q_{s,H}H,
\]
with $G_s=1$ on a zero-dimensional space. The right side of \eqref{eq:explicit-factors} is jointly measurable: the determinants and inverses may be computed on the fixed space $H$ using the positive definite operator $I_H-Q_{s,H}+Q_{s,H}\Sigma Q_{s,H}$. Empty products equal one. In the reverse case, the hypotheses on the input imply in addition that
\begin{equation}\label{eq:explicit-positive}
 h_j>0\quad\gamma_{A_j}\text{-almost everywhere whenever }\quad
 \mu\{s:p(s)<0,\ j\in J_s\}>0.
\end{equation}
For a family of the form \eqref{eq:explicit-factors}, condition \eqref{eq:explicit-positive} is also sufficient for the negative-exponent marginals to be positive and finite almost everywhere. It does not replace the other hypotheses of the main theorems.
Conversely, data satisfying \eqref{eq:explicit-algebra}--\eqref{eq:explicit-factors} and the hypotheses of the relevant inequality yield equality. The corresponding compatible density is
\begin{equation}\label{eq:explicit-joint}
 r(z)=G(z_H)\prod_j h_j(z_j),
\end{equation}
where $G$ is the density of $N(b,\Sigma)$ with respect to $\gamma_H$.
\end{theorem}

\begin{proof}
We first prove necessity. Equality in either main theorem gives a compatible density $r$. Put $r_t=\mathcal Q_t r$. The compatibility independence properties hold on a common set of full measure for all positive times, with the marginals smoothed by $\mathcal Q_t$. For each $s$ in that set, independence across $Q_s\R^n$ and its orthogonal complement gives
\begin{equation}\label{eq:explicit-commute}
 Q_s\nabla^2\log r_t(z)=\nabla^2\log r_t(z)Q_s,
 \qquad \nabla\log r_t(z)\in V.
\end{equation}
Here the commutation is on $\R^n$, and both identities hold for every $z$. Intersect the full-measure set from compatibility with $S_0$. Its projections still span $\mathcal L$: otherwise, their span would be a proper subspace containing $Q_s$ for almost every $s$, contrary to the definition of $\mathcal L$. Choose a finite spanning family from this intersection and denote the selected projections by $Q_1',\ldots,Q_d'$. Their complete common eigenspaces are precisely $F_1,\ldots,F_\ell$.

We use the argument in Lemma~\ref{lem:support}. Fix $z$ and write $C_v=\partial_v(\nabla^2\log r_t)(z)$. Symmetry of third derivatives and differentiation of the commutation identities give
\[
 C_vw=C_wv,\qquad Q_i'C_v=C_vQ_i',\qquad
 C_{Q_i'v}=Q_i'C_v.
\]
Since $C_{Q_i'v}$ also commutes with $Q_j'$, it follows that
\[
 (Q_i'Q_j'-Q_j'Q_i')C_v=0.
\]
The space $W=\operatorname{span}\{C_vw:v,w\in\R^n\}$ is invariant under every $Q_i'$, and their restrictions to $W$ commute. They therefore have a common orthonormal eigenbasis on $W$. Because the eigenvalues of projections are zero or one, $W\subset F_1\oplus\cdots\oplus F_\ell$. Each symmetric $C_v$ has range in $W$, and so vanishes on $D=(F_1\oplus\cdots\oplus F_\ell)^\perp$. Thus the $D$-block of the Hessian is constant. The Hessian preserves every $F_j$, and by symmetry also $D$. Hence the mixed blocks vanish. Integrating these identities yields
\begin{equation}\label{eq:explicit-smooth}
 \log r_t(z)=c_t+\tfrac12\langle D_tz_D,z_D\rangle
 +\langle d_t,z_D\rangle+\sum_{j=1}^\ell\psi_{j,t}(P_{F_j}z),
\end{equation}
where $D=(F_1\oplus\cdots\oplus F_\ell)^\perp$ and $D_t$ is symmetric on $D$. In particular, $\nabla^2\log r_t$ is block diagonal for $D,F_1,\ldots,F_\ell$, and its $D$-block is constant.

For $u,v\in F_j$, keeping the other components fixed in \eqref{eq:explicit-smooth} gives
\[
 \nabla\psi_{j,t}(u)-\nabla\psi_{j,t}(v)\in F_j\cap V=A_j.
\]
It follows that the projection of $\nabla\psi_{j,t}$ onto $F_j\cap A_j^\perp$ is constant. Thus $\psi_{j,t}$ is the sum of a function of $P_{A_j}z$ and a linear function of the orthogonal component in $F_j$. Consequently $r_t\gamma_n$ is a product of a Gaussian probability measure on $H$ and probability measures on the spaces $A_j$. Let $b_t$ and $\Sigma_t$ be the mean and covariance of the Gaussian factor on $H$. Since $r_t$ is a positive smooth probability density, $\Sigma_t$ is positive definite.

The subspace $A$ is contained in $V$, and therefore $P_HV\subset V_H$. Projecting $\nabla\log r_t(z)\in V$ onto $H$, and using the Gaussian density on $H$, shows that
\[
 (I_H-\Sigma_t^{-1})H\subset V_H,\qquad
 \Sigma_t^{-1}b_t\in V_H.
\]
The first operator is symmetric and has range in $V_H$, so it vanishes on $H\cap V_H^\perp$. Hence $\Sigma_t^{-1}$ is the identity on $H\cap V_H^\perp$ and leaves $V_H$ invariant. The same holds for $\Sigma_t$. Consequently $(\Sigma_t-I_H)H\subset V_H$, and $\Sigma_t^{-1}b_t\in V_H$ implies $b_t\in V_H$. Since $H$ reduces every $Q_s$ in the chosen full-measure set, restricting \eqref{eq:explicit-commute} to $H$ gives $Q_{s,H}(I_H-\Sigma_t^{-1})=(I_H-\Sigma_t^{-1})Q_{s,H}$. Thus $\Sigma_tQ_{s,H}=Q_{s,H}\Sigma_t$ on the same set.

Fix one $t>0$ and write $a=e^{-t}$, $\sigma^2=1-a^2$. If $Y\sim r\gamma_n$ and $Z\sim\gamma_n$ is independent, then $aY+\sigma Z$ has law $r_t\gamma_n$. For $\xi=\xi_H+\sum_j\xi_j$ in the orthogonal decomposition, independence at time $t$ gives
\[
 \mathbb E e^{ia\langle\xi,Y\rangle}e^{-\sigma^2|\xi|^2/2}
 =\mathbb E e^{i\langle\xi_H,aY_H+\sigma Z_H\rangle}
 \prod_j\mathbb E e^{i\langle\xi_j,aY_j+\sigma Z_j\rangle}.
\]
Substituting each marginal characteristic function and cancelling the nonzero Gaussian factors yields
\[
 \mathbb E e^{ia\langle\xi,Y\rangle}
 =\mathbb E e^{ia\langle\xi_H,Y_H\rangle}
 \prod_j\mathbb E e^{ia\langle\xi_j,Y_j\rangle}.
\]
Since $a>0$, uniqueness of characteristic functions proves that $Y_H,Y_1,\ldots,Y_\ell$ are independent. The $H$-component has characteristic function of Gaussian form with
\[
 b=a^{-1}b_t,\qquad \Sigma=a^{-2}(\Sigma_t-\sigma^2I_H).
\]
If $v$ is an eigenvector of $\Sigma$ with a negative eigenvalue, the modulus of the purported characteristic function along $\lambda v$ exceeds one, a contradiction. Hence $\Sigma\ge0$. If a nonzero $v$ belongs to its kernel, the characteristic function of $\langle v,Y_H-b\rangle$ is identically one, so this random variable is zero almost surely. This would put the $H$-marginal on a proper affine hyperplane. Such concentration is impossible because $r\gamma_n$ and its $H$-marginal are absolutely continuous. Thus $\Sigma>0$ whenever $H\ne\{0\}$.  Moreover,
\[
 \Sigma-I_H=a^{-2}(\Sigma_t-I_H),
\]
so $b\in V_H$ and \eqref{eq:explicit-algebra} holds. Denote the densities of the $A_j$-marginals of $r\gamma_n$ by $h_j$. This proves \eqref{eq:explicit-joint}.

For almost every $s$, the projection $Q_s$ acts as zero or identity on each $A_j$, and it commutes with $\Sigma$ on $H$. Taking the $B_s$-marginal of \eqref{eq:explicit-joint} gives \eqref{eq:explicit-factors}. In the reverse case, if the set in \eqref{eq:explicit-positive} has positive measure, choose a negative-exponent parameter in this set for which \eqref{eq:explicit-factors} and the input positivity hold. If $h_j$ vanished on a set of positive $\gamma_{A_j}$-measure, the right side of \eqref{eq:explicit-factors} would vanish on a set of positive $\gamma_n$-measure, contradicting positivity of $g_s(B_sz)$. This proves necessity.

Conversely, suppose \eqref{eq:explicit-positive} holds. There are only finitely many $j$, so, after removing the union of the null sets $\{s:p(s)<0,\ j\in J_s\}$ corresponding to factors not required to be positive, every factor seen by a remaining negative-exponent parameter is positive almost everywhere. The factors are also finite almost everywhere because they are probability densities. Fubini's theorem and positivity of $G_s$ then give $0<g_s<\infty$ $\gamma_k$-almost everywhere for almost every negative-exponent parameter. This proves the stated sufficiency of \eqref{eq:explicit-positive} for the positivity requirement.

Conversely, suppose the stated data exist and define $r$ by \eqref{eq:explicit-joint}. It is a probability density. Since $b\in V_H$ and $\Sigma-I_H$ is symmetric with range in $V_H$, it vanishes on $H\cap V^\perp$ and leaves $V_H$ invariant. The corresponding Gaussian cross-covariance is zero. Hence the $H$-Gaussian has an independent standard Gaussian component on $H\cap V^\perp$. The independent $A_j$ factors lie in $V$. Hence the $V$ and $V^\perp$ components have the required compatibility property. For almost every $s$, each $A_j$ lies either in the range or in the kernel of $Q_s$, while \eqref{eq:explicit-algebra} gives
\[
 Q_{s,H}\Sigma(I_H-Q_{s,H})=0.
\]
The two projected components of the $H$-Gaussian are therefore independent. Together with the independent $A_j$ factors, this proves that $Q_sY$ and $(I-Q_s)Y$ are independent for $Y\sim r\gamma_n$. Finally, \eqref{eq:explicit-factors} identifies the $B_s$-marginal with $g_s\gamma_k$: $B_s$ is an isometry from $Q_s\R^n$ onto $\R^k$, and sends the relevant standard Gaussian reference measures to the corresponding standard Gaussian measures. Thus $r$ is compatible, and equality follows from Theorem~\ref{thm:forward} or Theorem~\ref{thm:reverse}.
\end{proof}

\begin{remark}
The arbitrary factors in Theorem~\ref{thm:explicit} occur only on the subspaces $A_j$. If all $A_j$ vanish, every equality case is Gaussian. The Gaussian covariance need not be scalar: it is determined by the linear constraints \eqref{eq:explicit-algebra}. The theorem gives a structure theorem within the finite-dimensional common-source and coisometry hypotheses of this paper. 
\end{remark}

\section{Quadratic functions}

We conclude with an application of the semigroup identity \eqref{eq:dissipation} to exponential quadratic functions in the mixed-exponent reverse regime, without requiring a uniform positive lower bound on the fiber precision matrices. Consider
\[
f_s(x) = \exp\left(c_s + \langle a_s, x \rangle + \frac{1}{2}\langle H_s x, x \rangle\right),
\qquad H_s = H_s^*, \quad J_s = I - p(s) H_s > 0,
\]
and assume $C \ge P$ together with the integrability condition
\begin{equation}\label{eq:quadratic-integrability}
\int_\Omega \left(|c_s| + |a_s|^2 + \|H_s\|^2 + \tr|\log J_s| + \langle J_s^{-1} a_s, a_s \rangle\right) d\mu(s) < \infty.
\end{equation}
Set $a = \int_\Omega B_s^* a_s\,d\mu(s)$ and $H = \int_\Omega B_s^* H_s B_s\,d\mu(s)$, and suppose that $I - H > 0$.

\begin{corollary}\label{cor:quadratic}
Under the assumptions above, the reverse inequality holds:
\[
\int_{\R^n} F_-(z) \exp\left(\int_\Omega \frac{\log m_s}{p(s)}\,d\mu(s)\right) d\gamma_n(z)
\ge \exp\left(\int_\Omega \frac{\log m_s}{p(s)}\,d\mu(s)\right).
\]
Equality holds if and only if
\begin{equation}\label{eq:quadratic-eq}
a \in V, \quad H(\R^n) \subset V, \qquad p(s) a_s = B_s a, \quad p(s) H_s B_s = B_s H \quad \text{for } \mu\text{-a.e. } s.
\end{equation}
\end{corollary}

\begin{proof}
For $t > 0$, write $\tau = 1 - e^{-2t} \in (0,1)$ and $D_s(t) = (1 - \tau)I + \tau J_s = I - \tau p(s) H_s > 0$. Completing the square in the Gaussian integral defining $\mathcal{Q}_t(f_s^{p(s)})$ yields
\[
u_s(t,x) = \frac{1}{p(s)} \log \mathcal{Q}_t(f_s^{p(s)})(x)
= c_s(t) + \langle a_s(t), x \rangle + \frac{1}{2}\langle H_s(t) x, x \rangle,
\]
where
\begin{align*}
H_s(t) &= e^{-2t} H_s D_s(t)^{-1}, \qquad a_s(t) = e^{-t} D_s(t)^{-1} a_s, \\
c_s(t) &= c_s - \frac{\log \det D_s(t)}{2 p(s)} + \frac{p(s) \tau}{2} \langle D_s(t)^{-1} a_s, a_s \rangle.
\end{align*}
For any finite $t > 0$, we have $D_s(t) \ge e^{-2t} I$. Using \eqref{eq:quadratic-integrability} together with $0 < D_s(t)^{-1} \le e^{2t} I$ and $|\tr \log D_s(t)| \le \tr|\log J_s|$ (which follows from comparing the eigenvalues $1 - \tau + \tau \lambda_j$ of $D_s(t)$ with $\lambda_j$ of $J_s$), the parameter integrals defining $U_t(z) = \int_\Omega u_s(t, B_s z)\,d\mu(s)$ and its derivatives are finite and admit differentiation under the integral sign.

Let $H(t) = \int_\Omega B_s^* H_s(t) B_s\,d\mu(s)$ and $a(t) = \int_\Omega B_s^* a_s(t)\,d\mu(s)$. At $t = 0$, $H(0) = H$ satisfies $I - H(0) > 0$. As long as $I - H(t) > 0$, the Gaussian integral $e^{\Phi(t)} = \int_{\R^n} e^{U_t}\,d\gamma_n$ is finite and equals
\[
\Phi(t) = \int_\Omega c_s(t)\,d\mu(s) - \frac{1}{2}\log\det(I - H(t)) + \frac{1}{2}\langle (I - H(t))^{-1} a(t), a(t) \rangle,
\]
and the dissipation identity \eqref{eq:dissipation} with $h_{t,z}(s) = a_s(t) + H_s(t) B_s z$ gives $\Phi'(t) \le 0$. If $I - H(t)$ were to become singular at a finite first time $t_* > 0$, then $-\frac{1}{2}\log\det(I - H(t)) \to +\infty$ as $t \uparrow t_*$ while $\int_\Omega c_s(t)\,d\mu(s)$ remains bounded and $\langle (I - H(t))^{-1} a(t), a(t) \rangle \ge 0$, contradicting $\Phi(t) \le \Phi(0)$. Hence $I - H(t) > 0$ for all $t \ge 0$, and $\Phi(t)$ is nonincreasing on $[0,\infty)$.

As $t \to \infty$ ($\tau \uparrow 1$), we have $D_s(t) \to J_s$, so $c_s(t) \to p(s)^{-1} \log \int_{\R^k} f_s^{p(s)}\,d\gamma_k$. Because $|\tr \log D_s(t)| \le \tr|\log J_s|$ and $\tau D_s(t)^{-1} \le J_s^{-1}$ (as $D_s(t) \ge \tau J_s$), dominated convergence gives $\int_\Omega c_s(t)\,d\mu(s) \to \int_\Omega p(s)^{-1}\log m_s\,d\mu(s)$. For the linear term, the arithmetic--geometric mean inequality gives $(1 - \tau)I + \tau J_s \ge 2\sqrt{(1-\tau)\tau}\,J_s^{1/2}$, whence
\[
|a_s(t)| = e^{-t} |D_s(t)^{-1/2} D_s(t)^{-1/2} a_s| \le \frac{1}{2\sqrt{\tau}} \langle J_s^{-1} a_s, a_s \rangle^{1/2},
\]
so $a(t) \to 0$ as $t \to \infty$ by dominated convergence. For the quadratic term, the eigenvalues of $p(s) H_s$ lie in $(-\infty, 1)$, so {for an eigenvalue $\lambda<1$ of $p(s)H_s$,
\[
 \left|\frac{(1-\tau)\lambda}{1-\tau\lambda}\right|
 \le \max\left\{1,\frac{1-\tau}{\tau}\right\}.
\]
Consequently $\|H_s(t)\|\le p_*^{-1}\max\{1,(1-\tau)/\tau\}$. For each $s$, $H_s(t)\to0$ as $t\to\infty$. Dominated convergence therefore gives $H(t)\to0$.} Thus $\lim_{t \to \infty} \Phi(t) = \int_\Omega p(s)^{-1}\log m_s\,d\mu(s)$, proving the inequality.

{If equality holds, the nonincreasing function $\Phi$ has equal values at $0$ and $\infty$, so $\Phi'(0)=0$. By \eqref{eq:dissipation} and \eqref{eq:kernel},
\[
 h_z:s\longmapsto a_s+H_sB_sz\in\Ker(P-C)
\]
for almost every $z\in\R^n$. The map $z\mapsto h_z$ is continuous from $\R^n$ into $\Hil$, and the kernel is closed. Hence the identity holds for every $z$. By \eqref{eq:kernel},
\[
 h_z(s)=\frac{1}{p(s)}B_sv(z),\qquad v(z)\in V.
\]
Applying $A^*$ gives $v(z)=A^*h_z=a+Hz$. Thus
\[
 a+Hz\in V,\qquad
 p(s)(a_s+H_sB_sz)=B_s(a+Hz)
\]
for every $z$ and almost every $s$, with a common exceptional set obtained by taking a countable dense set of $z$. Taking $z=0$ and subtracting the resulting identity proves \eqref{eq:quadratic-eq}.}

Conversely, suppose \eqref{eq:quadratic-eq} holds. Multiplying $p(s) H_s B_s = B_s H$ on the right by $B_s^*$ and using $B_s B_s^* = I_k$ gives $p(s) H_s = B_s H B_s^*$, and therefore $Q_s H = B_s^* B_s H = p(s) B_s^* H_s B_s = H Q_s$. Since $H(\R^n) \subset V$, we may restrict to an orthonormal eigenbasis of $H$ in $V$. On $V$, $K = \int_\Omega p(s)^{-1} Q_s\,d\mu(s) = I_V$. {The spectral theorem for the commuting symmetric operators $Q_s$ and $H$ gives
\[
 \log(I-Q_sH)=Q_s\log(I-H),\qquad
 (I-Q_sH)^{-1}Q_s=(I-H)^{-1}Q_s.
\]
These identities use only $I-H>0$ and do not require $\|H\|<1$. Hence}
\begin{align*}
	\int_\Omega \frac{\log \det J_s}{p(s)}\,d\mu(s)
&= \int_\Omega \frac{\tr \log(I - p(s) H_s)}{p(s)}\,d\mu(s)
= \int_\Omega \frac{\tr \log(I - Q_s H)}{p(s)}\,d\mu(s)
\\&= \tr\bigl(K \log(I - H)\bigr) = \log\det(I - H).
\end{align*}
{Similarly, using the second identity above and} $p(s) a_s = B_s a$ with $a \in V$ gives
\begin{align*}
	\int_\Omega \frac{\langle J_s^{-1} p(s) a_s, p(s) a_s \rangle}{p(s)}\,d\mu(s)
&= \int_\Omega \frac{\langle (I - Q_s H)^{-1} Q_s a, a \rangle}{p(s)}\,d\mu(s)
\\&= \langle (I - H)^{-1} K a, a \rangle = \langle (I - H)^{-1} a, a \rangle.
\end{align*}
Hence $\Phi(0) = \Phi(\infty)$, so equality holds.
\end{proof}

\begin{example}[Absence of a uniform precision bound]
Let $n = k = 1$ and $B_s = 1$ for all $s$. Take a positive atom $\{0\}$ of unit mass with $p(0) = 1/2$, and negative atoms $j \ge 1$ of mass $2^{-j}$ with $p(j) = -1$. Set $H_0 = 0$, $H_j = -1 + 1/j$ for $j \ge 1$, and $a_s = c_s = 0$. Then $J_j = 1 - (-1)(-1 + 1/j) = 1/j$, so the fiber precision matrices $J_j$ have no uniform positive lower bound. Nevertheless, $C \ge P$ and \eqref{eq:quadratic-integrability} hold, and Corollary~\ref{cor:quadratic} yields
\[
\int_\R \exp\left(\frac{1}{2}\sum_{j=1}^\infty 2^{-j}\left(-1 + \frac{1}{j}\right)x^2\right) d\gamma_1(x)
= (2 - \log 2)^{-1/2}
> \exp\left(-\frac{1}{2}\sum_{j=1}^\infty 2^{-j} \log j\right).
\]
Here $V = \R$ and $H = \sum_{j=1}^\infty 2^{-j}(-1 + 1/j) = \log 2 - 1$, so the equality condition $p(s) H_s = H$ fails and the inequality is strict.
\end{example}

\section*{\textbf{Acknowledgments}}
The authors developed the main ideas
and approach of this work and are responsible for the mathematical
results and proofs. ChatGPT 6 was used as an assistive tool in
exploring some proof strategies, checking calculations, and improving
the language and presentation of the manuscript.

\end{document}